\documentclass[11pt]{article}

\usepackage[a4paper,margin=2.5cm]{geometry}
\usepackage[T1]{fontenc}
\usepackage[utf8]{inputenc}
\usepackage{lmodern}
\usepackage{amsmath,amssymb,bm,mathtools}
\usepackage{graphicx}
\usepackage[detect-all]{siunitx}
\usepackage[export]{adjustbox}
\usepackage{subfig}
\usepackage{comment}
\usepackage{enumitem}
\usepackage{etoolbox}
\usepackage{hyperref}
\usepackage{xparse}

\ProvideDocumentCommand{\bmsection}{s m}{%
  \IfBooleanTF{#1}{\section*{#2}}{\section{#2}}%
}

\ProvideDocumentCommand{\bmsubsection}{s m}{%
  \IfBooleanTF{#1}{\subsection*{#2}}{\subsection{#2}}%
}

\DeclarePairedDelimiter{\norm}{\lVert}{\rVert}
\newcommand{\normtwo}[1]{\norm{#1}_{L_2}} % \normtwo{x} -> ||x||_2
\makeatletter
\patchcmd{\@makecaption}{\small}{\normalsize}{}{}
\makeatother

\title{A Quasicontinuum Method with Optimized Local Maximum-Entropy Interpolation and Heaviside Enrichment for Heterogeneous Lattices}
\author{Benjamin Werner$^{1,2}$\thanks{Corresponding author: \href{mailto:b.werner@tue.nl}{b.werner@tue.nl}} \and Ond\v{r}ej Roko\v{s}$^{2}$ \and Jan Zeman$^{1}$}
\date{}

\begin{document}
\maketitle

\begin{center}
$^{1}$Faculty of Civil Engineering, Czech Technical University in Prague, Th\'{a}kurova 7, 166 29 Prague 6, Czech Republic\\
$^{2}$Mechanics of Materials, Department of Mechanical Engineering, Eindhoven University of Technology, P.O. Box 513, 5600 MB Eindhoven, The Netherlands
\end{center}

\begin{abstract}
Lattice systems are indispensable for modeling and analyzing physical phenomena in materials with discrete or heterogeneous micro- or meso-structures. However, the computational requirements for practical engineering applications of lattice systems remain high. The QuasiContinuum (QC) method addresses this by reducing the system of equations using a finite element mesh interpolation, rather than considering all nodes of a fully resolved lattice. Nevertheless, interfaces between separate phases in heterogeneous materials, such as concrete, require fine meshes throughout the domain, diminishing the effectiveness of QC. Enrichment strategies originally introduced in the eXtended Finite Element Method (XFEM) can also account for material interfaces in discrete systems using nonconforming meshes, thereby resolving this issue. In a previous study, we applied Heaviside enrichment to investigate concrete mesostructures with the extended QC method, achieving a tenfold reduction in the number of unknowns while maintaining similar accuracy compared to discretizations with fully resolved interfaces. In the present study, we employ the meshless Local Maximum Entropy (LME) interpolation, which transitions seamlessly from widespread meshfree to linear basis functions. Additionally, we combine LME interpolation with Heaviside enrichment and systematically investigate the role of the locality parameter and its optimization in heterogeneous lattices. This combination of optimized LME basis functions with Heaviside enrichment leads to an order-of-magnitude improvement in displacement accuracy while using the same number of Degrees of Freedom (DOF) compared to QC with linear interpolation. Moreover, we identify optimized distributions of the LME locality parameter and propose simple, non-optimized rules that deliver comparable accuracy at a fraction of the computational cost. Results from three numerical examples show that the optimal locality-parameter fields are nonuniform near interfaces and can be approximated by simple pattern-based rules that retain much of the benefit of full optimization.
\end{abstract}

\noindent\textbf{Keywords:} QuasiContinuum Methodology, Local Maximum-Entropy Interpolation, Optimum Support Size, Weak Discontinuities, XFEM, Discrete Lattice Systems

\section{Introduction}\label{sec1}

%\subsection{Lattice models and their application}
Discrete lattice systems assembled from truss or beam elements are essential for the computational modeling of various materials. Examples include textiles and woven fabrics \cite{CHENG2016170,BEEX201382,BENBOUBAKER2007498}, in which the elements represent the yarns and the lattice nodes represent the inter-yarn junctions. They have been used to study the failure of fibrous materials such as paper \cite{LIU20101} as well as cortical bovine bone \cite{Mayya2016}. More recently, Braun et al. \cite{BRAUN2021107767,BRAUN2024105213} used discrete lattice networks to study damage and fracture in fiber-reinforced composites. In addition, discrete lattice systems have a long history in numerical modeling of failure in quasi-brittle materials such as concrete and rock \cite{Nikolic2018753,NIKOLIC2015209}. They are relatively easy to implement and numerically robust for simulating failure phenomena, such as crack initiation, growth, and branching, at the meso- and micro-scales. The crack propagation is a natural outcome of the breakage of the lattice trusses or beam elements, which is beneficial to study quasi-brittle heterogeneous materials and materials with complex crack patterns. Furthermore, they overcome the disadvantages of the Finite Element Method (FEM), such as stress locking, strength dependence of the element size, and fracture toughness \cite{BOLANDER2021108030,LilliuPhD2007}. 

%\subsection{QC method}
However, a major disadvantage of discrete lattice systems is their high computational cost. Because of the small and fixed lattice spacing, they do not allow coarsening at locations of low interest, which differs substantially compared to the numerical models of continua. This enforces a similar size of lattice members throughout the entire numerical model and leads to a large number of degrees of freedom (DOF). The multiscale QuasiContinuum (QC) method has proven to be effective in reducing the computational cost of discrete lattice systems and was pioneered by Tadmor et al. to study dislocations and nanoindentation in atomistic lattices \cite{Tadmor19961529,Tadmor19964529}. The method reduces computational cost through two key components: \textit{interpolation} and \textit{summation}. The QC method reduces the number of unknowns by interpolating lattice node (atom) displacements over a coarse finite-element mesh rather than resolving all lattice DOFs. Summation rules further reduce the cost by selecting a subset of sampling atoms or interactions with associated weights to approximate the total potential energy of the system. The weights indicate the number of sites that the sampling atom or interaction represents in an element \cite{AMELANG2015378}. With these two integral parts, the QC method enables the study of solids at the atomistic level at locations with highly non-uniform deformations, while seamlessly transitioning to the continuum level at places with small deformation gradients, thus limiting computational efforts. \cite{Kochmann2016}

%\subsection{Applications and extensions of the QC method}
The QC method was originally developed to couple atomistic and continuum models efficiently and has since been extended to a wide range of materials and mechanical phenomena. It has been applied on the microscale to study the friction behavior of rough aluminum surfaces \cite{TRAN2019180} and, at the mesoscopic and macroscopic levels, to study fabrics with elastic-plastic material behavior.\cite{BEEX201552} The QC method holds potential for irregular lattice arrangements, which is of great advantage for use in rock and concrete analysis \cite{MIKES201750}. It has been extended to beam lattice structures and referred to as ``generalized QC,'' with nonlinear elastic~\cite{PHLIPOT2019758} and plastic~\cite{CHEN2020112878,Chen20212498} properties to predict their deformation behavior, as well as to study crack initiation and fracture toughness~\cite{Kraschewski20241432}. QC has been extensively used to predict fracture in atomistic models of metals.\cite{Qiu2017} A generalization to dissipative processes such as fracture based on a variational formulation of QC for spring lattices was introduced by Roko\v{s} et al.\cite{Rokos2017174,ROKOS2016214} 

%\subsection{XFEM for strong and weak discontinuities}
Enriched methods of the eXtended Finite Element Method (XFEM) greatly simplify the discretization of structures and materials with strong and weak discontinuities. \cite{Moes2003,Fries2010,Cheng2010,Li2018,Zhu2012186,WU201577,GUPTA201323} Cracks or material interfaces in heterogeneous materials do not have to be resolved by the FE mesh and instead are considered by additional DOFs and enriched interpolation functions. This allows crack paths that are unknown a priori to be predicted and eliminates mesh dependence in fracture mechanics problems of continua. In QC, these enriched methods have been utilized to also coarsen the mesh in the wake of a growing crack to avoid a trail of a fully-resolved region and reduce the computational burden tremendously for fractured spring lattice structures. \cite{ROKOS2017769} In a recent study, we have applied enrichment strategies for weak discontinuities together with a nonconforming mesh to account for material interfaces. \cite{Werner2024e7415} This allows for the determination of homogenized material properties of highly heterogeneous discrete systems with accuracy similar to that of a conforming mesh, with low computational costs.

%\subsection{Meshless and maximum entropy interpolation}
A variety of meshless methods have been developed over the last decades that avoid the use of elements for approximation, in contrast to the Finite Element Method (FEM). They overcome strong connectivity constraints of the FEM and nodes can be placed more freely in a domain. This simplifies resolving material interfaces and crack tips, and avoids a remeshing for such problems. \cite{AlMahmoud2024725}
The Local Maximum Entropy (LME) approximation scheme is one such meshless method that is characterized by an optimal compromise between two competing interests: (i) the least biased estimate (maximum entropy principle) by having wide stretching interpolation functions (i.e., large support) and (ii) local shape functions of least width (i.e., small support).\cite{Arroyo20062167} The LME shape functions have global maximum-entropy and affine interpolation functions as limiting cases. The locality parameter of the functions allows a smooth transition between both of these extreme cases. Moreover, the LME interpolation satisfies the weak Kronecker-delta property at convex boundaries and thus facilitates the application of displacement boundary conditions. Arroyo and Ortiz \cite{Arroyo20062167} showed that the accuracy of the LME interpolation scheme is vastly superior to that of finite element shape functions.

LME interpolation has been applied to solve phase-field models to predict fracture in Kirchhoff-Love thin shells \cite{Amiri2014102} as well as to solve the Cahn-Hilliard equation for phase separation \cite{Amiri20191}. In addition, they have been used to study fracture in three-dimensional single edge notch beam models \cite{Amiri2016254} and lead to great agreement compared to experimental results. Millan et al. \cite{Millan2011723} used maximum entropy interpolation for the Kirchhoff-Love thin shell analysis to predict deformations of smooth manifolds under static loading. The interpolation scheme provides great flexibility in investigating surfaces of complex topology, while exhibiting better accuracy than finite elements. Foca \cite{Foca2015} used LME functions for thermo-mechanical simulations of rotary friction welding. The meshless interpolation is able to predict the large deformations during the welding process. Kumar et al. \cite{Kumar2019858} address tensile stability, which is one of the major drawbacks of multiple meshfree interpolation methods including the LME scheme. Tensile instability occurs from changing nodal spacing under finite deformations and is associated with the localization of the interpolation function support. They present a stable maximum entropy interpolation scheme by using an updated Lagrangian formulation and an anisotropic adaptation method for the interpolation functions based on the deformation gradient. This provides an accurate and efficient way to avoid tensile instability under large deformations. Rosolen et al. \cite{Rosolen2010868} improved the approximation accuracy of LME further by optimizing the locality parameter at each node. They show highly accurate solutions with very coarse discretizations for various problems. In addition, the optimizations reveal unexpected patterns of the shape functions support size. Kochmann and Venturini \cite{Kochmann2014034007} applied the maximum-entropy interpolation method to the QC method with efficient quadrature type summation rules. Due to the affine interpolation as a limiting case of maximum entropy functions, the approximation scheme is able to recover the full atomistic resolution and therefore allows a seamless transition between interpolated domains and fully resolved regions. The suggested adaptation concept \cite{Kochmann2014034007} thereby overcomes a limitation of previously used mesh free interpolation schemes in the QC context. Amiri et al. \cite{Amiri201445} used the LME approximation scheme in combination with enrichment functions from the XFEM to study an infinite plate with center crack under tensile loading. The approximation scheme predicts stress intensity factors with one or more orders of magnitude greater accuracy than XFEM. Rosolen and Arroyo \cite{Rosolen201395} combined the LME interpolation scheme with isogeometric analysis. The two interpolation schemes complement each other by eliminating their shortcomings. Using an isogeometric representation along the boundary leads to an exact representation of the geometry and facilitates the application of essential boundary conditions on non-convex boundaries. Applying maximum entropy interpolation in the bulk of the domain preserves the flexibility of using unstructured grids of points and the seamless transition between affine and widespread interpolation functions. They illustrate the capabilities of the method for heat equation and nonlinear elastic problems. Fathi et al. \cite{Fathi20216103} extended isogeometric analysis (IGA) using LME functions for the enriched interpolation (X-IGALME). The method is applied to two-dimensional fracture mechanics problems, where IGA has the advantage of being able to approximate complex boundaries. Using LME for the enhanced part leads to an improved estimation of stresses and a better crack propagation direction compared to XFEM or XIGA. 

%subsection{Gap}
Even with the high accuracy reported and the promising applications demonstrated, several open questions remain for LME. For heterogeneous materials, there is no systematic study of how to choose the spatial distribution of the LME locality parameter, and thus the shape-function support, near inclusions and interfaces, nor of simple rules that avoid costly optimization for each node. When LME interpolation is combined with enrichment to represent weak discontinuities, the selection or optimization of the locality parameter has not been assessed for heterogeneous materials. The integration of enriched LME interpolation into heterogeneous lattice networks within a QC framework remains essentially unexplored.

%subsection{Outlook}
%To close this gap, in this study, we employ the LME interpolation within the QC framework to systematically investigate its performance and the effect of the locality parameter within heterogeneous lattice networks. We combine LME shape functions with Heaviside enrichment to represent material interfaces while retaining a regular grid of repatoms, thereby reducing the number of unknowns. The LME locality parameter $\boldsymbol{\gamma}^{\mathrm{LME}}$, which controls the support of the basis function, is optimized following Rosolen et al.\ \cite{Rosolen2010868}. This strategy yields roughly an order of magnitude fewer degrees of freedom than discretizations using affine FE interpolation, while achieving comparable displacement accuracy. The present work focuses on interpolation; the development and analysis of a suitable summation rule are deliberately left for future studies.

To close this gap, in this study, we employ the LME interpolation within the QC framework to systematically investigate its performance and the effect of the locality parameter within heterogeneous lattice networks. We combine LME shape functions with Heaviside enrichment to represent material interfaces while retaining a regular grid of repatoms, thereby reducing the number of unknowns. A central objective is to show that the optimal locality-parameter fields are nonuniform near interfaces and that their structure can be characterized systematically. Based on these observations, we identify simple pattern-based rules that recover much of the benefit of full optimization. This provides a practical alternative to expensive per-repatom optimization. The LME locality parameter $\boldsymbol{\gamma}^{\mathrm{LME}}$, which controls the support of the basis function, is optimized following Rosolen et al.\ \cite{Rosolen2010868}. This strategy yields roughly an order of magnitude fewer degrees of freedom than discretizations using affine FE interpolation, while achieving comparable displacement accuracy. The present work focuses on interpolation; the development and analysis of a suitable summation rule are deliberately deferred to future studies.

\section{Methods}\label{sec2}

\subsection{Full-lattice formulation}

\begin{figure}[ht]
	\centering
	\begin{minipage}{0.25\linewidth}
		\begin{center}
			\includegraphics[scale=1]{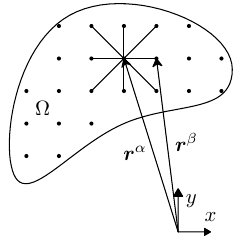}
		\end{center}
		\vspace{-1em}
		global kinematic variable: \\
			$\mathbf{r}=[\mathbf{r}^1,\dots,\mathbf{r}^{n_\mathrm{ato}}]^\mathsf{T}$
	\end{minipage}
	\hspace{0.1\linewidth}
	\begin{minipage}{0.4\linewidth}	
		\begin{center}
			\includegraphics[scale=1]{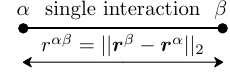}
		\end{center}
		\vspace{-0.9em}
		elastic strain: $e^{\alpha\beta}(r^{\alpha\beta}) = (r^{\alpha\beta}-r_0^{\alpha\beta})/r_0^{\alpha\beta}$\\
	\end{minipage}
	\caption{\normalsize Schematic illustration of the kinematic variables of an X-braced lattice.}
	\label{fig:Schematic_Kinematic}
\end{figure}

In this study, two-dimensional X-braced lattice arrangements are considered. Each lattice node, which we will refer to as an atom, has eight direct neighbors connected by interactions within the domain $\Omega$ (Figure \ref{fig:Schematic_Kinematic}). 
The location of the atoms is stored in the position vector as a single column matrix $\mathbf{r} = [\mathbf{r}^1,\ldots,\mathbf{r}^{n_\mathrm{ato}}]^T$, $\mathbf{r}^{\alpha} \in \mathbb{R}^2, \alpha = 1, \ldots, n_{\mathrm{ato}}$ for the entire system with $n_\mathrm{ato}$ atoms. A pseudo-time parameter $t_k$ ($k=0,\ldots,n_{\mathrm{T}}$) is introduced, leading to a discretization of the time interval $[0,T]$ as $0 = t_0 < t_1 < \cdots < t_{n_{\mathrm{T}}} = T$ in a quasi-static setting. By minimizing the total potential energy $\Pi^k$ at a pseudo-time instant $t_k$ 
\begin{equation}
    \mathbf{r}(t_{k}) = \underset{\widehat{\mathbf{r}} \in \mathbb{R}^{2\, n_{\mathrm{ato}}}}{\mathrm{arg\,min\,}} \Pi^k(\widehat{\mathbf{r}}), \quad k=1,\ldots,n_{\mathrm{T}},
\end{equation}
an equilibrium configuration $\mathbf{r}$ of the lattice system is determined, using the initial condition $\mathbf{r}(0) = \mathbf{r}_0$. Here, $\widehat{\mathbf{r}}$ denotes an arbitrary admissible position vector. The total potential energy of the lattice system
\begin{equation}
    \Pi^k(\widehat{\mathbf{r}}) = \sum \limits_{\alpha \beta = 1}^{n_{\mathrm{int}}} \pi^k_{\alpha\beta} (\widehat{\mathbf{r}})
\end{equation}
is the sum of the elastic energies $\pi^k_{\alpha \beta}$ of all interactions at time $t_k$. The lattice network has $n_{\mathrm{int}}$ interactions, stored in an index set $N_{\mathrm{int}}$, and the interaction $\alpha\beta \in N_{\mathrm{int}}$ connects atoms $\alpha$ and $\beta$. The length of the interaction $\alpha\beta$
\begin{equation}
    r^{\alpha\beta} = \normtwo{\mathbf{r}^{\beta} - \mathbf{r}^{\alpha}},
\end{equation} 
is defined by the Euclidean norm of the difference between the two position vectors $\mathbf{r}^{\alpha}$ and $\mathbf{r}^{\beta}$. The elastic energy of the interaction between atoms $\alpha$ and $\beta$ is given by
\begin{equation}
    \pi^k_{\alpha\beta} (\widehat{\mathbf{r}}) = \dfrac{E^{\alpha\beta}A^{\alpha\beta}}{2r_0^{\alpha\beta}}(\widehat{r}^{\alpha\beta}-r_0^{\alpha\beta})^2
    \label{eq:Interaction_ElEn}
\end{equation}
where $E^{\alpha\beta}$ is the Young's modulus and $A^{\alpha\beta}$ is the cross-sectional area of the interaction. In addition, the strain of the interaction, as the difference between the initial $r^{\alpha\beta}_0$ and the current length $\widehat{r}^{\alpha\beta}$, divided by the initial length, is part of the interaction elastic energy in Eq. (\ref{eq:Interaction_ElEn}). The governing equilibrium equations follow from the stationarity of the total potential energy,
\begin{equation}
    \dfrac{\partial \Pi(\widehat{\mathbf r})}{\partial \widehat{\mathbf r}} = \mathbf{0}
\end{equation}
and are solved iteratively by the Newton–Raphson method. The equilibrium of energy stationarity and the Newton–Raphson solution strategy; for more detail, see Werner et al.~\cite{Werner2024e7415}.

%%%%%%%%%%%%%%%%%%%%%%%%%%%%%%%%%%%%%%%%%%%%%%%%%%%%%%%%%%%%%%%%%%%%%%%%%%%%%
%% QC reduction by interpolation
%%%%%%%%%%%%%%%%%%%%%%%%%%%%%%%%%%%%%%%%%%%%%%%%%%%%%%%%%%%%%%%%%%%%%%%%%%%%%
\subsection{QC reduction using LME interpolation and Heaviside enrichment}

In a discrete lattice system, the QC method reduces the number of unknowns by interpolating atomic displacements usually with linear \cite{Tadmor19961529,Tadmor19964529,Rokos2017174,ROKOS2017769} or higher-order finite elements \cite{Beex2015344,Kraschewski20241432}. In the present study, we employ LME shape functions enriched with Heaviside interpolation to represent interfaces (see Werner et al. \cite{Werner2024e7415} for more details). A set of representative atoms (repatoms) with positions $\widehat{\mathbf{r}}_{\mathrm{rep}}$ serve as interpolation nodes. The position vector of atom $\alpha$
\begin{equation}
     \widehat{\mathbf{r}}^{\alpha}_{\mathrm{qc}} 
     = 
     \sum_{\beta \in N_{\mathrm{rep}}} \phi_{\beta}(\mathbf{r}_0^{\alpha})\widehat{\mathbf{r}}_{\mathrm{rep}}^{\beta}
     +
    \sum_{j=1}^{n^{\star}} 
    \underbrace{
    \phi_{\beta_j}(\mathbf{r}_0^{\alpha})[\chi(\mathbf{r}_0^{\alpha}) - \chi(\mathbf{r}_0^{\beta_j})]
    }_{\phi_j^{\star}(\mathbf{r}_0^{\alpha})}
    \widehat{\mathbf{g}}_{j}^{\star}, \quad \alpha \in N_{\mathrm{ato}}, \quad \beta_j \in N^{\star}
    \label{eq:interpoaltion}
\end{equation}
is determined by the LME shape functions $\phi_{\beta}(\mathbf{r}_0^{\alpha})$ and the locations of the repatoms $\widehat{\mathbf{r}}_{\mathrm{rep}}^{\beta}$ with $N_{\mathrm{ato}}$ and $N_{\mathrm{rep}}$ being index sets that include all $n_{\mathrm{ato}}$ atoms and $n_{\mathrm{rep}}$ repatoms. The enriched interpolation is calculated using the shape function $\phi_{\beta_j}(\mathbf{r}_0^{\alpha})$, with $\beta_j$ being one of the $n^{\star}$ enriched repatoms in the index set $N^{\star}$. The enriched degrees of freedom $\widehat{\mathbf{g}}_{j}^{\star}$ together with the Heaviside function in its shifted form $\chi(\mathbf{r}_0^{\alpha}) - \chi(\mathbf{r}_0^{\beta_j})$ account for weak discontinuities. The shifted Heaviside function ensures the Kronecker-delta property on the the boundary and greatly simplifies the application of boundary conditions. Multiplying  the shifted Heaviside function by the LME interpolation function $\phi_{\beta_j}$ yields the enriched interpolation function $\phi_j^{\star}(\mathbf{r}_0^{\alpha})$. The LME basis functions and the enriched interpolation functions are combined to the interpolation matrix
\begin{equation}
    \boldsymbol{\Phi} = [\boldsymbol{\Phi}_{\mathrm{LME}}, \boldsymbol{\Phi}^{\star}],
    \label{eq:InterpolationMatrix}
\end{equation}
and both are assembled as
\begin{equation}
    (\Phi_{\mathrm{LME}})_{(2\alpha-1)(2j-1)} = 
    \begin{cases}
        \phi_{\beta_j}(\mathbf{r}_0^{\alpha}), & \mathrm{for}\; \alpha \in N_{\mathrm{ato}}, \beta_j \in N_{\mathrm{rep}}, j=1,\ldots,n_{\mathrm{rep}} \\
        0, & \mathrm{otherwise}
    \end{cases}
    \label{eq:Phi_assembling}
\end{equation}
\begin{equation}
    (\Phi^{\star})_{(2\alpha-1)(2j-1)} =
    \begin{cases}
        \phi_{\mathrm{GSON}, \, j}^{\star}(\mathbf{r}_0^{\alpha}), & \mathrm{for}\; \alpha \in N_{\mathrm{ato}}, j=1,\ldots,n^{\star} \\
        0, & \mathrm{otherwise}
    \end{cases}
\end{equation}
exhibiting a checkerboard pattern of zero and nonzero values. Subsections~\ref{subsec:LME} and \ref{subsec:enrichment} summarize the computation of the LME interpolation and its key properties, and present the derivation of the enriched functions.

%%%%%%%%%%%%%%%%%%%%%%%%%%%%%%%%%%%%%%%%%%%%%%%%%%%%%%%%%%%%%%%%%%%%%%%%%%%%%
%% Local maximum-entropy (LME) shape functions
%%%%%%%%%%%%%%%%%%%%%%%%%%%%%%%%%%%%%%%%%%%%%%%%%%%%%%%%%%%%%%%%%%%%%%%%%%%%%
\subsection{Local maximum-entropy (LME) shape functions}
\label{subsec:LME}
Among meshfree methods, LME interpolation is characterized by its adjustable locality parameter $\beta^{\mathrm{LME}}$, which controls the support size of the shape function. This flexibility includes the limiting case of linear interpolation and enables an accurate interface resolution. Moreover, its Kronecker-delta property on the boundary allows the direct application of displacement boundary conditions without additional considerations \cite{Arroyo20062167}. However, inside the interpolated domain, the LME basis lacks the Kronecker–delta property. The DOFs of the repatoms are not equal to the displacements of the co-located atoms in the deformed configuration. Therefore, repatoms should be viewed as generalized degrees of freedom or control variables that parameterize the displacement field. The LME interpolation parameters and basis functions are defined in the undeformed reference configuration, where $\mathbf{r}_{0,\mathrm{rep}}^{\beta}$ denotes the position of the repatom located at atom~$\beta$. For brevity, we will nevertheless refer to the repatom located at atom~$\beta$ simply as repatom~$\beta$ in the following.

The LME shape function $\phi_{\alpha}(\mathbf{r}_{0})$ for the repatom $\alpha$ is given by \cite{Rosolen2010868}
\begin{equation}
    \phi_{\alpha}(\mathbf{r}_{0}) 
    = 
    \dfrac{1}{Z\left(\mathbf{r}_{0},\boldsymbol{\lambda}^*, \boldsymbol{\beta}^{\mathrm{LME}}\right)}
    \mathrm{exp} \left[ -\beta^{\mathrm{LME}}_{\alpha} \normtwo{\mathbf{r}_{0} - \mathbf{r}_{0, \, \mathrm{rep}}^{\alpha}}^2 
    +
    \boldsymbol{\lambda}^*(\mathbf{r}_{0}) \cdot \left(\mathbf{r}_{0} - \mathbf{r}_{0, \, \mathrm{rep}}^{\alpha}\right) \right].
    \label{eq:LME_shapefunction}
\end{equation}
The locality parameter $\beta^{\mathrm{LME}}_{\alpha}$ of repatom $\alpha$ is scaled by the distances between each atom in $N_{\mathrm{ato}}$ and repatom $\alpha$ in the initial configuration. The second term in Eq. (\ref{eq:LME_shapefunction}) is the inner product between the Lagrange multiplier $\boldsymbol{\lambda}^*(\mathbf{r}_{0})$ and the spatial separation of each atom in $N_{\mathrm{ato}}$ and repatom $\alpha$. The denominator in Eq. (\ref{eq:LME_shapefunction}), called the partition function
\begin{equation}
    Z\left(\mathbf{r}_{0},\boldsymbol{\lambda}^*, \boldsymbol{\beta}^{\mathrm{LME}}\right)
    =
    \sum_{\beta \in N_{\mathrm{rep}}} \mathrm{exp} \left[ -\beta^{\mathrm{LME}}_{\beta} \normtwo{\mathbf{r}_{0} - \mathbf{r}_{0, \, \mathrm{rep}}^{\beta}}^2 
    +
    \boldsymbol{\lambda}^*(\mathbf{r}_{0}) \cdot \left(\mathbf{r}_{0} - \mathbf{r}_{0, \, \mathrm{rep}}^{\beta}\right) \right]
\end{equation}
is the sum of all repatoms, where the distances between the atoms in $N_{\mathrm{ato}}$ and the repatom $\beta$ are multiplied by the locality parameter $\beta^{\mathrm{LME}}_{\beta}$ of the repatom $\beta$.

The locality parameter can be defined in three different ways: (i) $\beta^{\mathrm{LME}}$ can be set to the same value for all repatoms, (ii) $\boldsymbol{\beta}^{\mathrm{LME}}$ can be defined separately for each repatom as $\boldsymbol{\beta}^{\mathrm{LME}} = [\beta^{\mathrm{LME}}_1,...,\beta^{\mathrm{LME}}_{n_{\mathrm{rep}}}]^{\mathrm{T}}$, leading to different isotropic shape functions as described above, or (iii) $\boldsymbol{\beta}^{\mathrm{LME}}$ can be defined as a vector for each spatial dimension at each repatom, i.e., $\boldsymbol{\beta}^{\mathrm{LME}} = [\boldsymbol{\beta}^{\mathrm{LME}}_1,...,\boldsymbol{\beta}^{\mathrm{LME}}_{n_{\mathrm{rep}}}]^{\mathrm{T}}$, which allows anisotropic shape functions with different support sizes in individual directions \cite{Kumar2019858}. In this study, we limit ourselves to options (i) and (ii). Furthermore, multiplying $\beta^{\mathrm{LME}}$ by the repatom spacing $h$, results in a dimensionless value $\gamma^{\mathrm{LME}}_{\alpha} = \beta^{\mathrm{LME}}_{\alpha} h^2$, and is comparable between different repatom spacings.

When computing the LME shape functions in Eq.~(\ref{eq:LME_shapefunction}), the Lagrange multiplier
\begin{equation}
  \boldsymbol{\lambda}^*(\mathbf{r}_{0})
  = \operatorname*{arg\,min}_{\boldsymbol{\lambda}\in\mathbb{R}^{2}}
    \log Z(\mathbf{r}_{0},\boldsymbol{\lambda},\boldsymbol{\beta}^{\mathrm{LME}})
  \label{eq:min_lambda}
\end{equation}
is determined by this unconstrained minimization problem, with $^*$ being the optimized value. In Appendix~\ref{appendix:app1} a solution procedure using a regularized Newton–Raphson scheme is described.

\begin{figure}[tbp]
   \centering
   \subfloat[center repatom]{\includegraphics[width=0.33\textwidth]{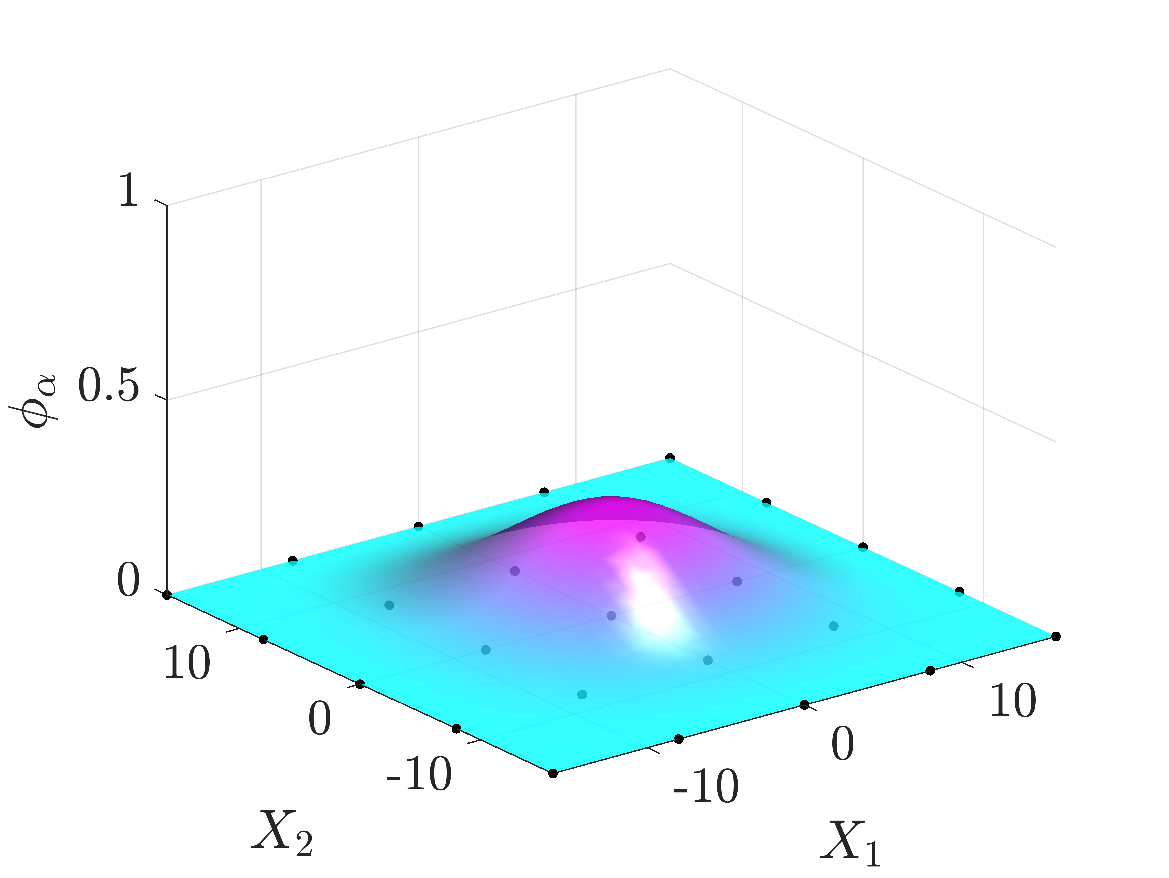}}
   \hfill
   \subfloat[boundary repatom]{\includegraphics[width=0.33\textwidth]{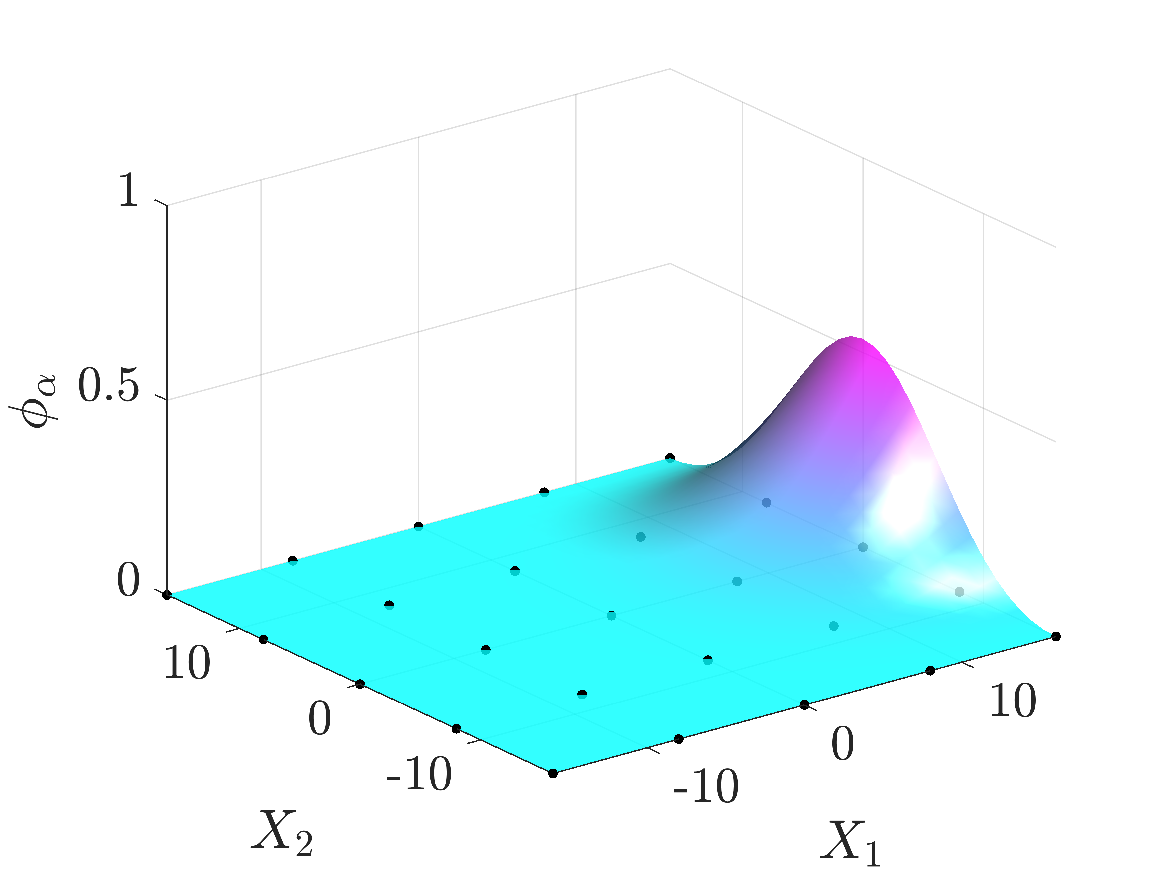}}
   \hfill
   \subfloat[corner repatom]{\includegraphics[width=0.33\textwidth]{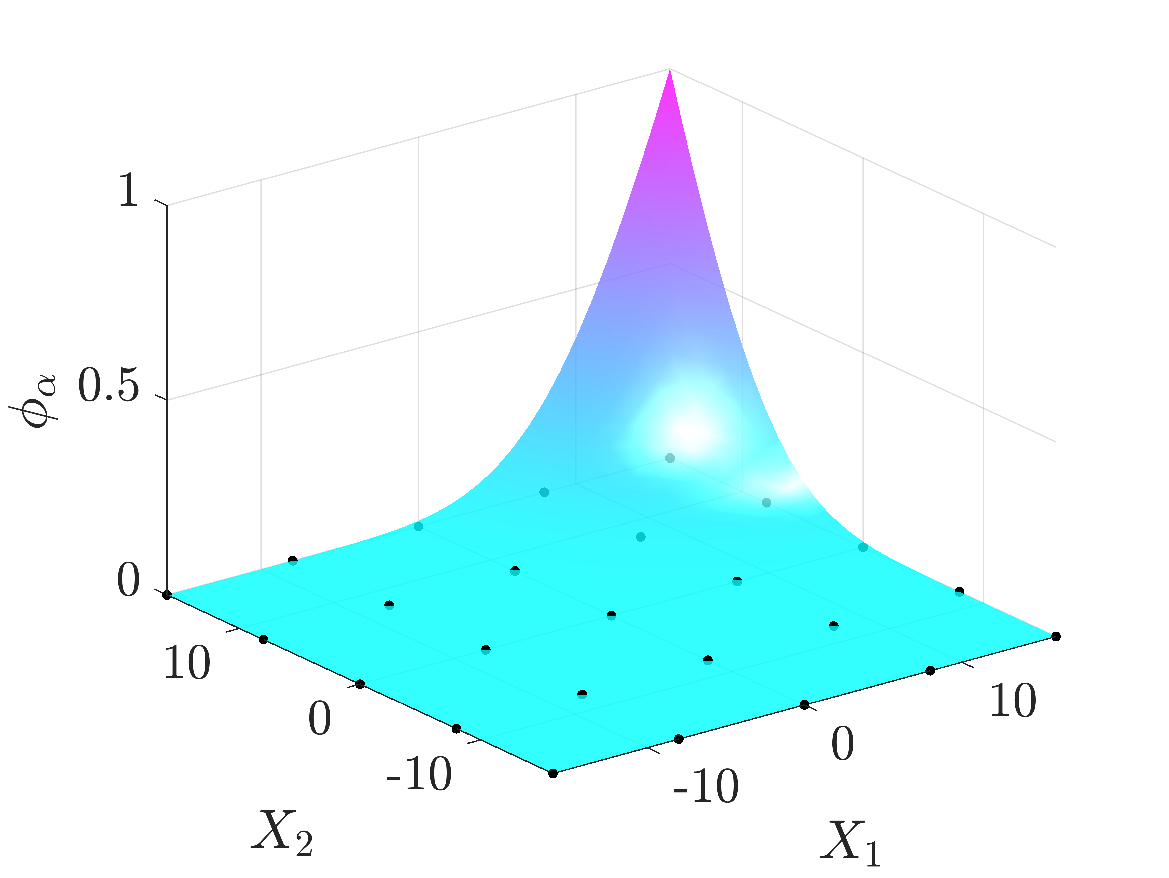}}
   \caption{Local Maximum Entropy interpolation functions with $\gamma^{\mathrm{LME}} = 0.9$ for  (a) the center repatom,  (b) a boundary repatom, and (c) a corner repatom of the two dimensional domain. The repatoms are indicated by black dots}
   \label{fig:Schematic_IntFunc_gamma09}
\end{figure}

\begin{figure}[tbp]
   \centering
   \subfloat[center repatom]{\includegraphics[width=0.33\textwidth]{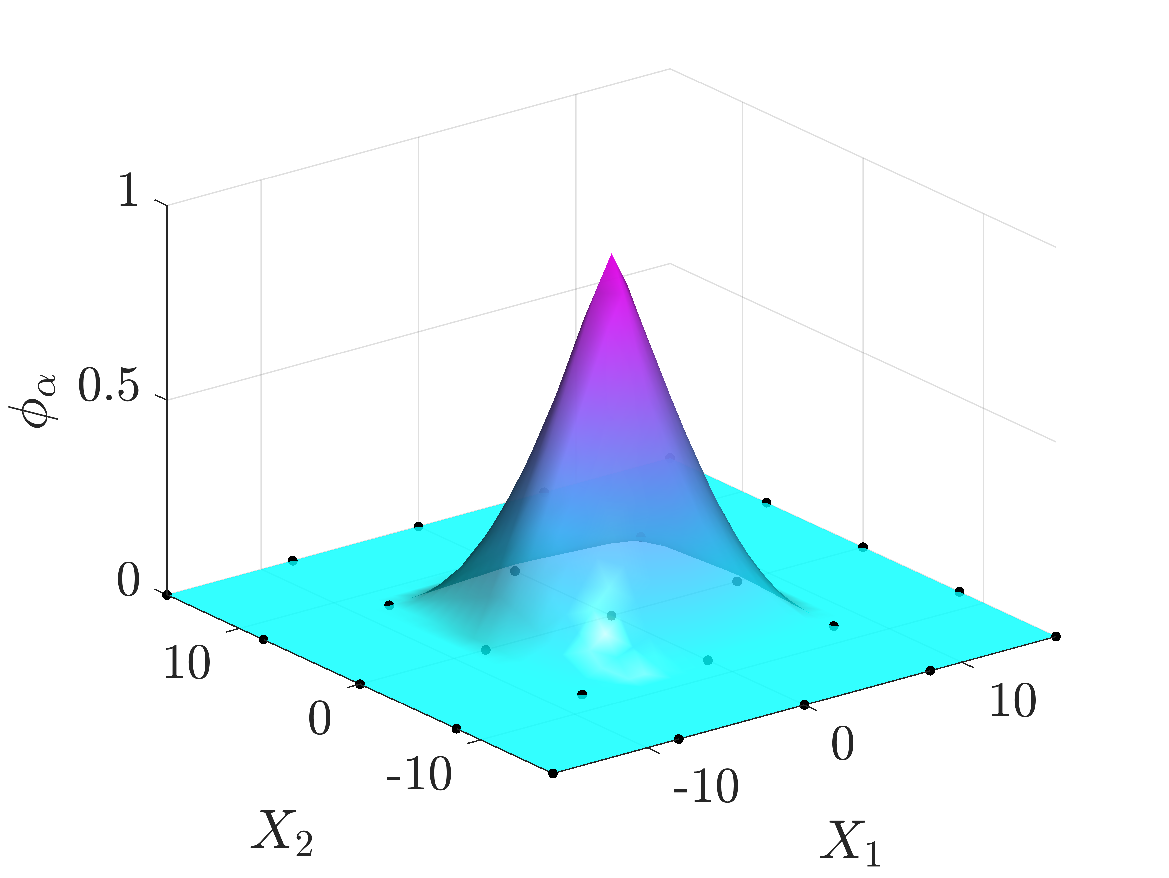}}
   \hfill
   \subfloat[boundary repatom]{\includegraphics[width=0.33\textwidth]{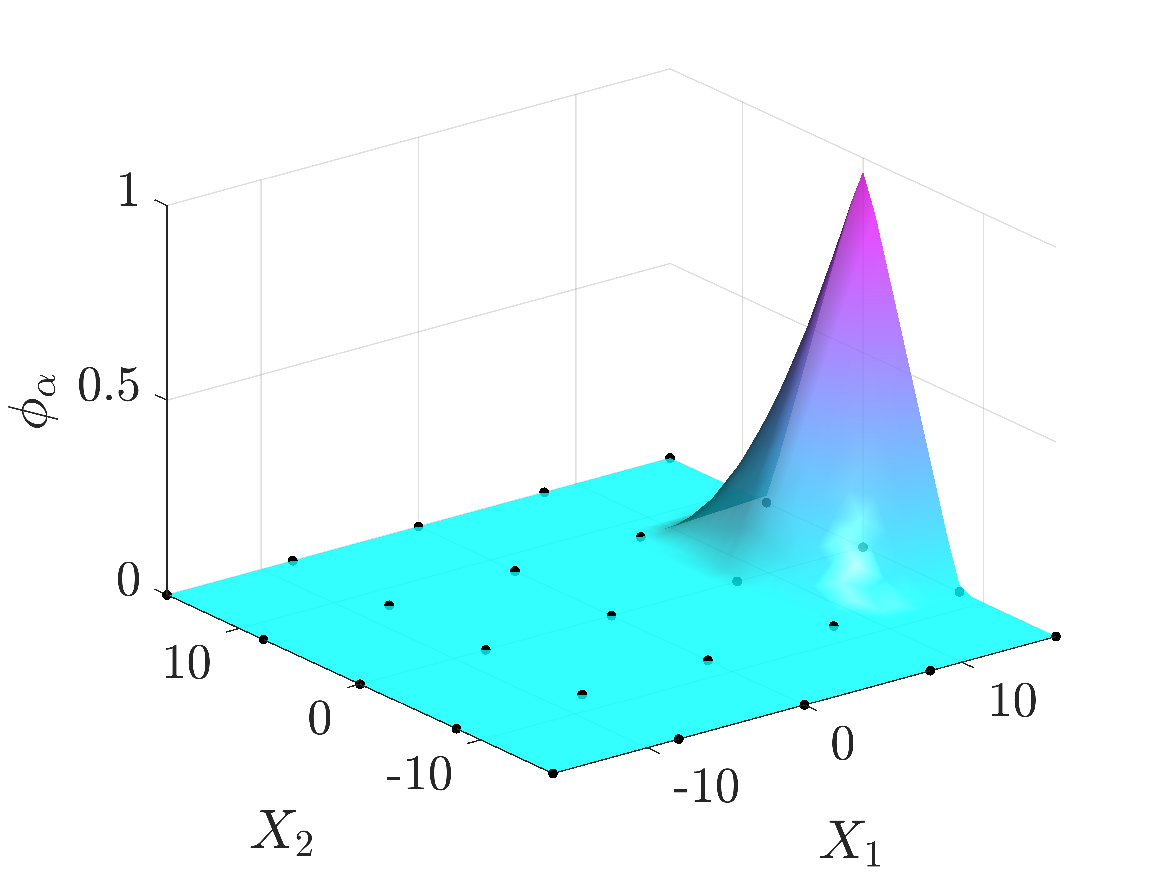}}
   \hfill
   \subfloat[corner repatom]{\includegraphics[width=0.33\textwidth]{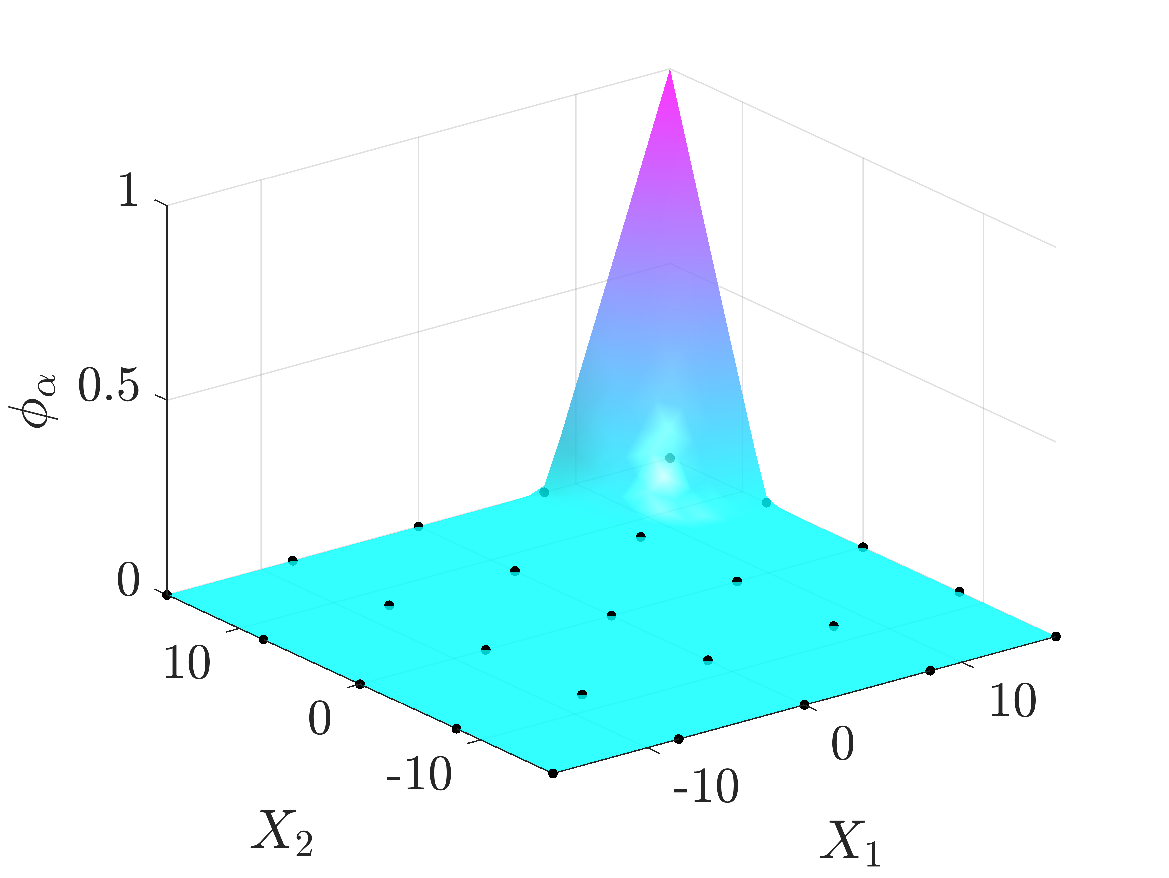}}
   \caption{Local Maximum Entropy interpolation functions with $\gamma^{\mathrm{LME}} = 4.0$ for  (a) the center repatom,  (b) a boundary repatom, and (c) a corner repatom of the two dimensional domain. The repatoms are indicated by black dots.}
   \label{fig:Schematic_IntFunc_gamma40}
\end{figure}

In Figures \ref{fig:Schematic_IntFunc_gamma09} and \ref{fig:Schematic_IntFunc_gamma40} the LME shape functions are schematically illustrated for a two-dimensional square domain, parameterized with 25 equally spaced repatoms. They are non-negative and satisfy the weak Kronecker delta property at the boundary of a convex hull of the domain. The interpolation functions of repatoms in the interior vanish at the boundary and the sum of the interpolation functions of repatoms along the boundary are approximately one, which simplifies applying displacement boundary conditions. Furthermore, in Figure \ref{fig:Schematic_IntFunc_gamma09} the basis functions have a large support due to the low dimensionless locality parameter $\gamma^{\mathrm{LME}}=0.9$, whereas the shape functions in Figure \ref{fig:Schematic_IntFunc_gamma40} with $\gamma^{\mathrm{LME}}=4.0$ are close to linear interpolation and have small support.

%%%%%%%%%%%%%%%%%%%%%%%%%%%%%%%%%%%%%%%%%%%%%%%%%%%%%%%%%%%%%%%%%%%%%%%%%%%%%
%% Heaviside enrichment and Gram--Schmidt orthogonalization
%%%%%%%%%%%%%%%%%%%%%%%%%%%%%%%%%%%%%%%%%%%%%%%%%%%%%%%%%%%%%%%%%%%%%%%%%%%%%
\subsection{Heaviside enrichment and Gram--Schmidt orthonormalization}
\label{subsec:enrichment} 
\begin{figure}[tbp]
   \centering
   \subfloat[Lattice with inclusion]{\includegraphics[width=0.27\textwidth]{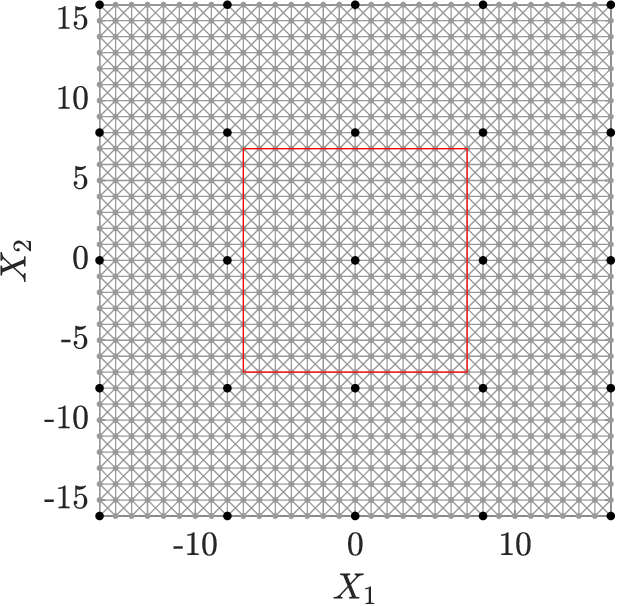}}
   \hfill
   \subfloat[Signed distance function]{\includegraphics[width=0.33\textwidth]{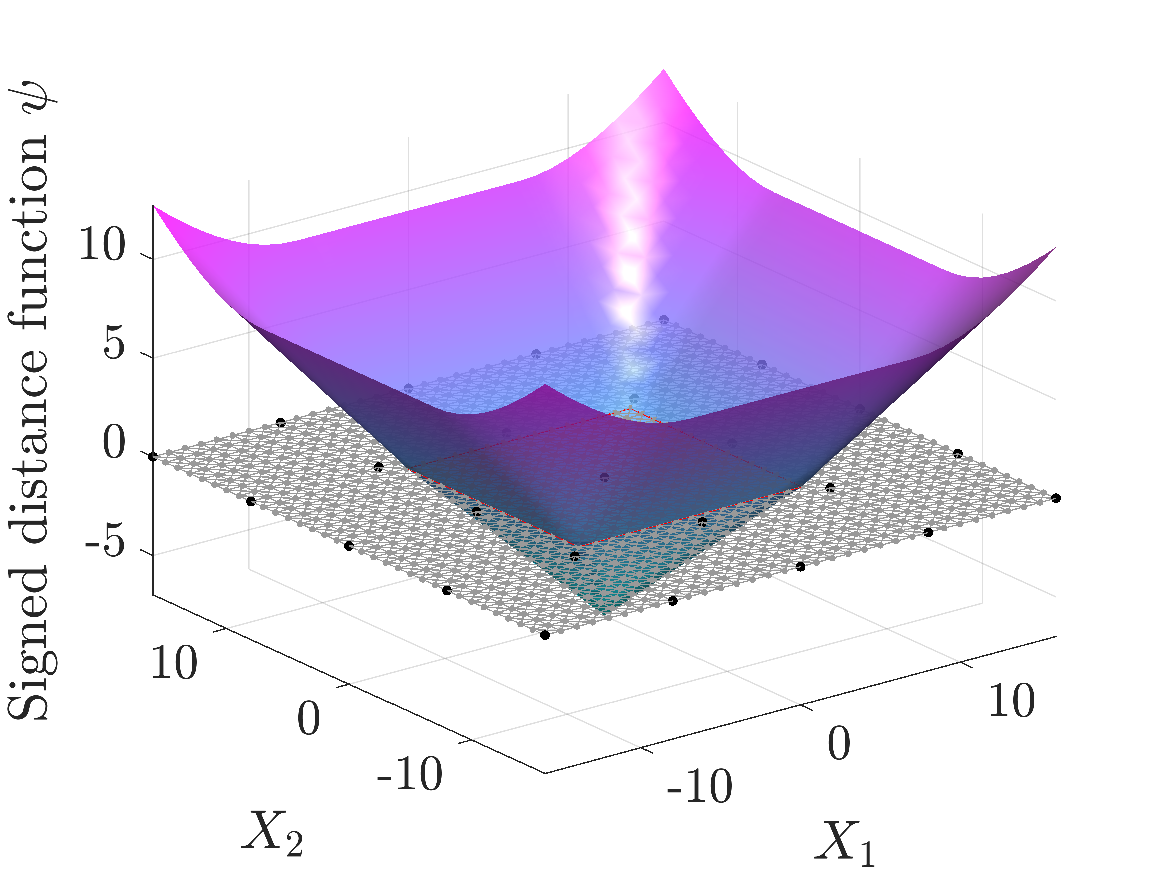}}
   \hfill
   \subfloat[Heaviside function]{\includegraphics[width=0.33\textwidth]{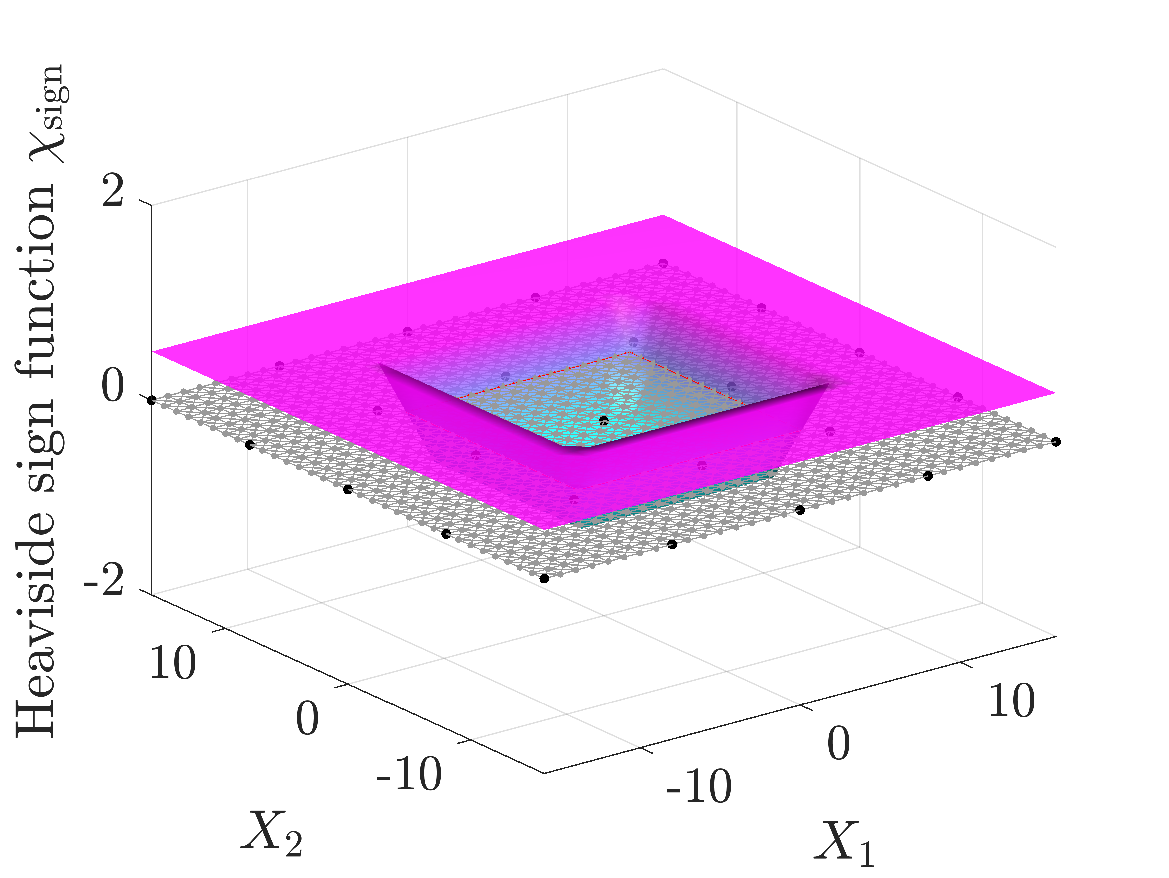}}
   \caption{(a) Square inclusion example with the interface shown in red, the underlying X-braced lattice in grey, and repatoms in black; (b) the corresponding signed distance function; (c) the Heaviside function derived from the signed distance function.}
   \label{fig:Schematic_SquareIncl}
\end{figure}

To construct the enriched interpolation functions, we extend the LME basis with a Heaviside function \cite{Moes1999} to account for weak discontinuities in discrete lattices. To illustrate construction, we use an example of a square inclusion example (Figure~\ref{fig:Schematic_SquareIncl}~a). The X-braced lattice is shown, with the inclusion interface highlighted in red, and the domain is parameterized by 25 repatoms, indicated by black dots in Figure~\ref{fig:Schematic_SquareIncl}~a.
The Heaviside sign function $\chi_{\mathrm{sign}}(\mathbf{r}_0^{\alpha})$ is determined from the signed distance function of the interface (Figure \ref{fig:Schematic_SquareIncl}~b)
\begin{equation}
    \psi(\mathbf{r}_0^{\alpha}) = \normtwo{\mathbf{r}_0^{\alpha} - \mathbf{c}^{\alpha}} \cdot \mathrm{sign}(\mathbf{n}_{c^{\alpha}}^{\mathrm{T}}(\mathbf{r}_0^{\alpha} - \mathbf{c}^{\alpha}))
    \label{eq:SignDistFunc}
\end{equation}
and calculated by the Euclidean norm between the location of the atom $\alpha$ in the initial configuration and the nearest location of an atom at the interface $\mathbf{c}^{\alpha}$. It is multiplied by the sign of the normal vector of the interface $\mathbf{n}_{c^{\alpha}}^{\mathrm{T}}$ of the atom $\alpha$ and is positive for atoms outside the inclusion and negative for all atoms inside the interface. The Heaviside sign function 
\begin{equation}
    \chi_{\mathrm{sign}}(\mathbf{r}_0^{\alpha}) = 
    \begin{cases}
        -0.5, & \psi(\mathbf{r}_0^{\alpha}) < 0, \\
        0, & \psi(\mathbf{r}_0^{\alpha}) = 0,\\
        +0.5, & \psi(\mathbf{r}_0^{\alpha}) > 0,
    \end{cases}
    \label{eq:enrichment_Heav}
\end{equation}
takes zero for every atom along the interface and is defined as $+0.5$ for atoms outside and $-0.5$ inside the inclusion (Figure \ref{fig:Schematic_SquareIncl}~c).

\begin{figure}[tbp]
   \centering
   \subfloat[center repatom]{\includegraphics[width=0.33\textwidth]{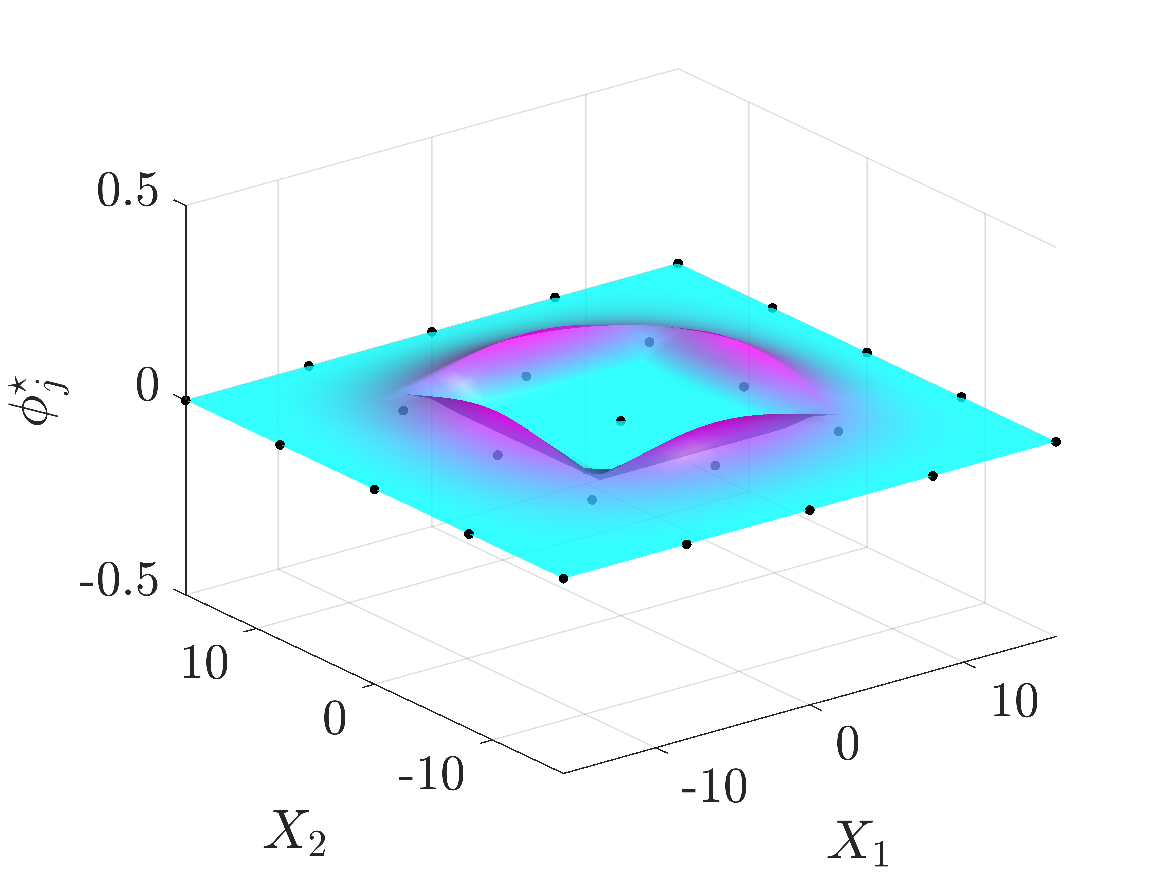}}
   \hfill
   \subfloat[boundary repatom]{\includegraphics[width=0.33\textwidth]{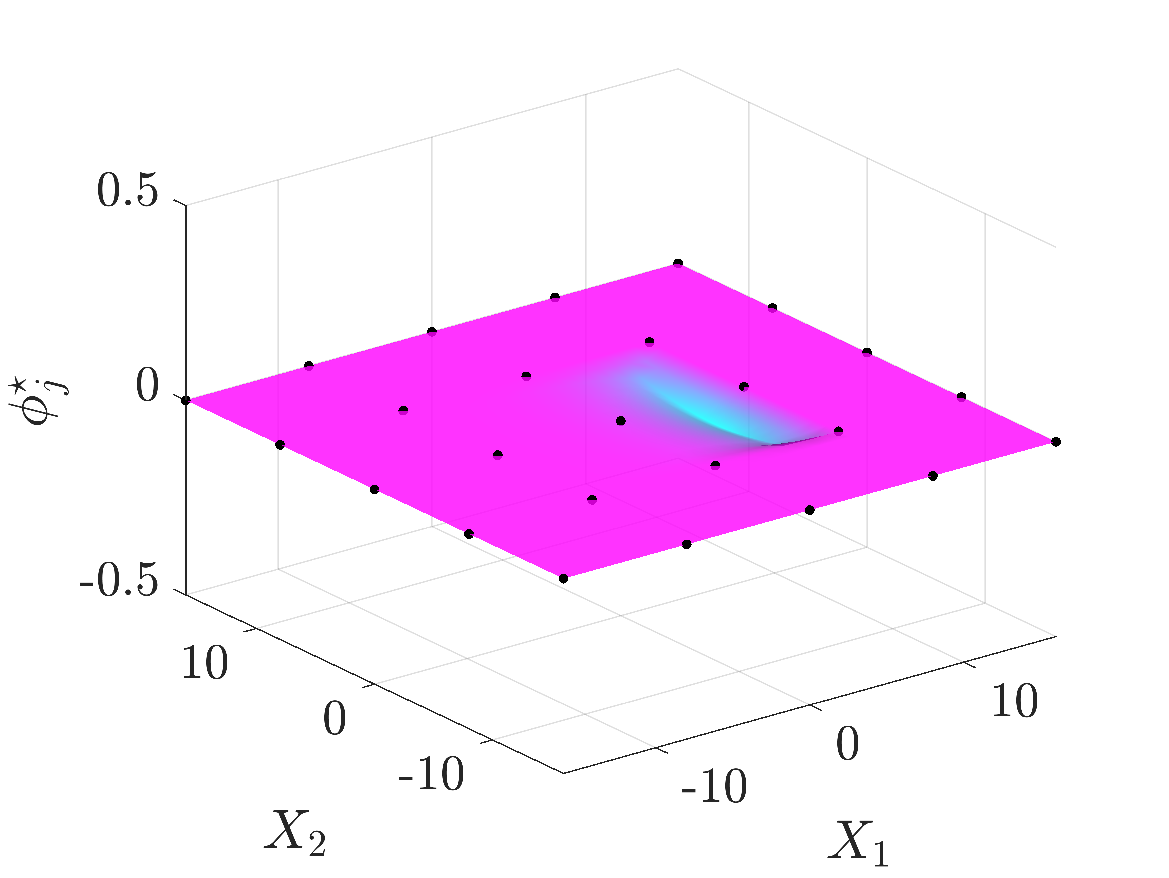}}
   \hfill
   \subfloat[corner repatom]{\includegraphics[width=0.33\textwidth]{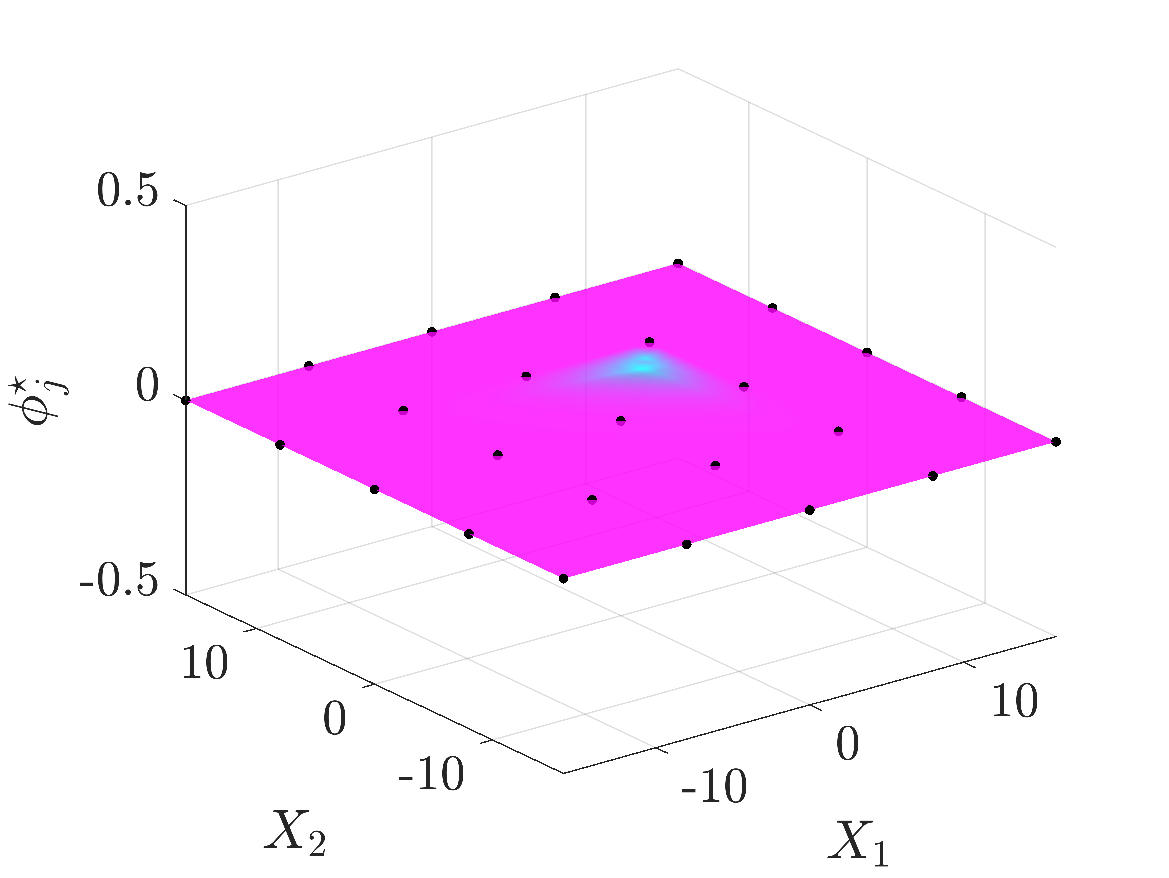}}
   \caption{Enriched LME interpolation functions using the shifted Heaviside sign function with $\gamma = 0.9$ for  (a) the center repatom,  (b) a boundary repatom, and (c) a corner repatom of the two dimensional domain. The repatoms are indicated by black dots}
   \label{fig:Schematic_IntFuncHeav_gamma09}
\end{figure}

\begin{figure}[tbp]
   \centering
   \subfloat[center repatom]{\includegraphics[width=0.33\textwidth]{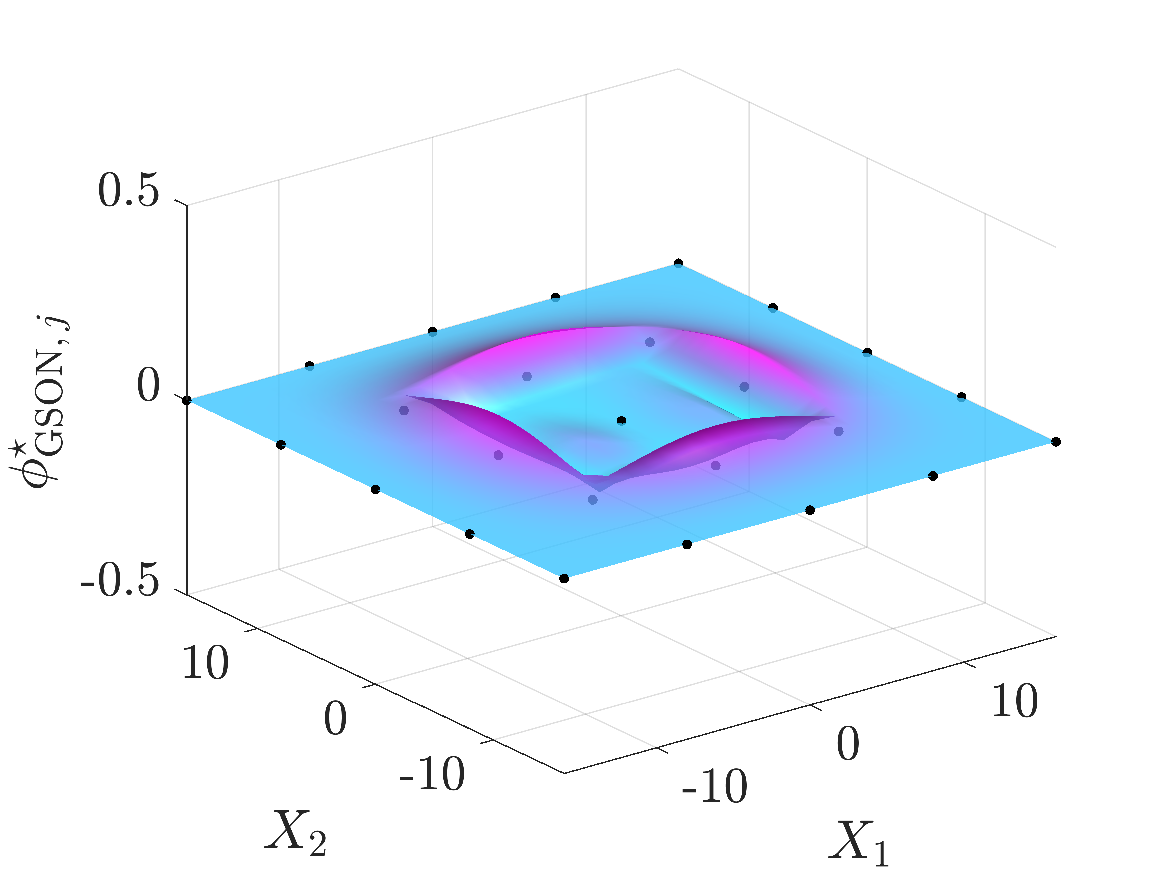}}
   \hfill
   \subfloat[boundary repatom]{\includegraphics[width=0.33\textwidth]{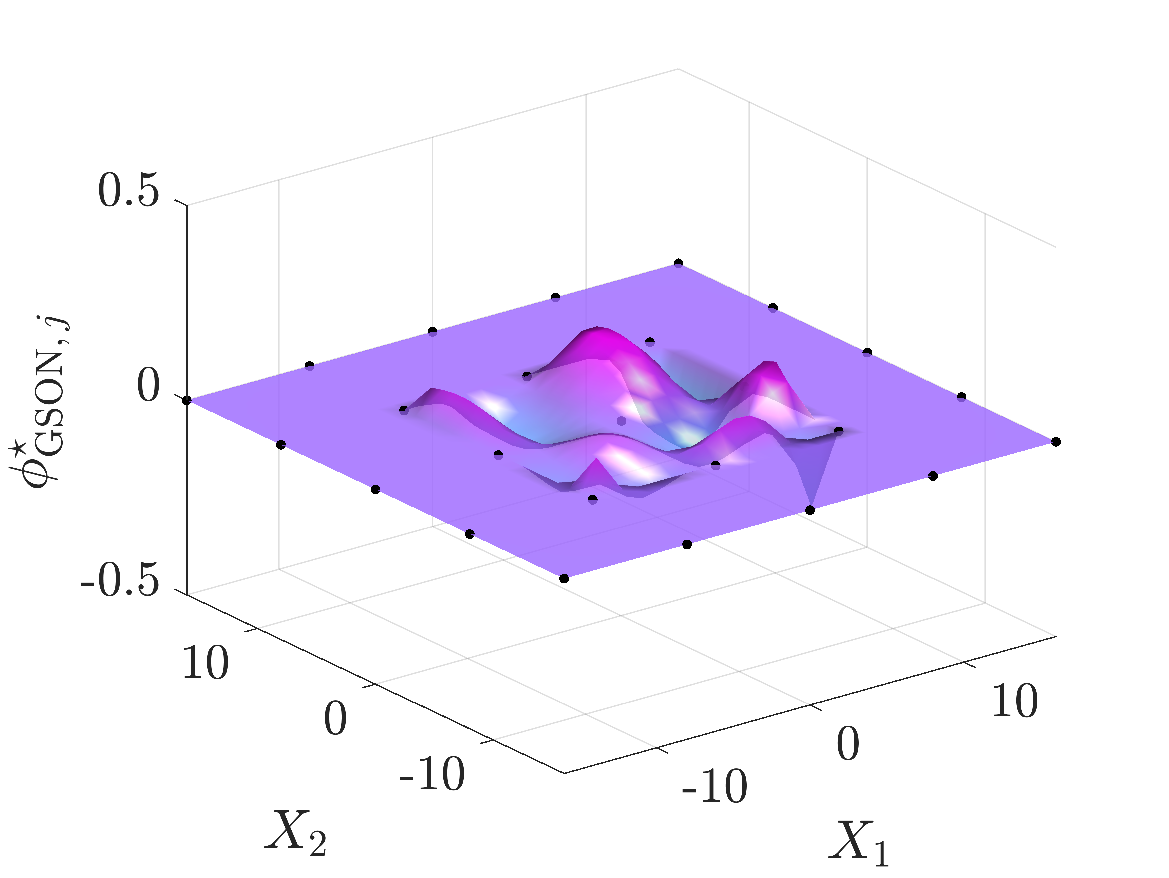}}
   \hfill
   \subfloat[corner repatom]{\includegraphics[width=0.33\textwidth]{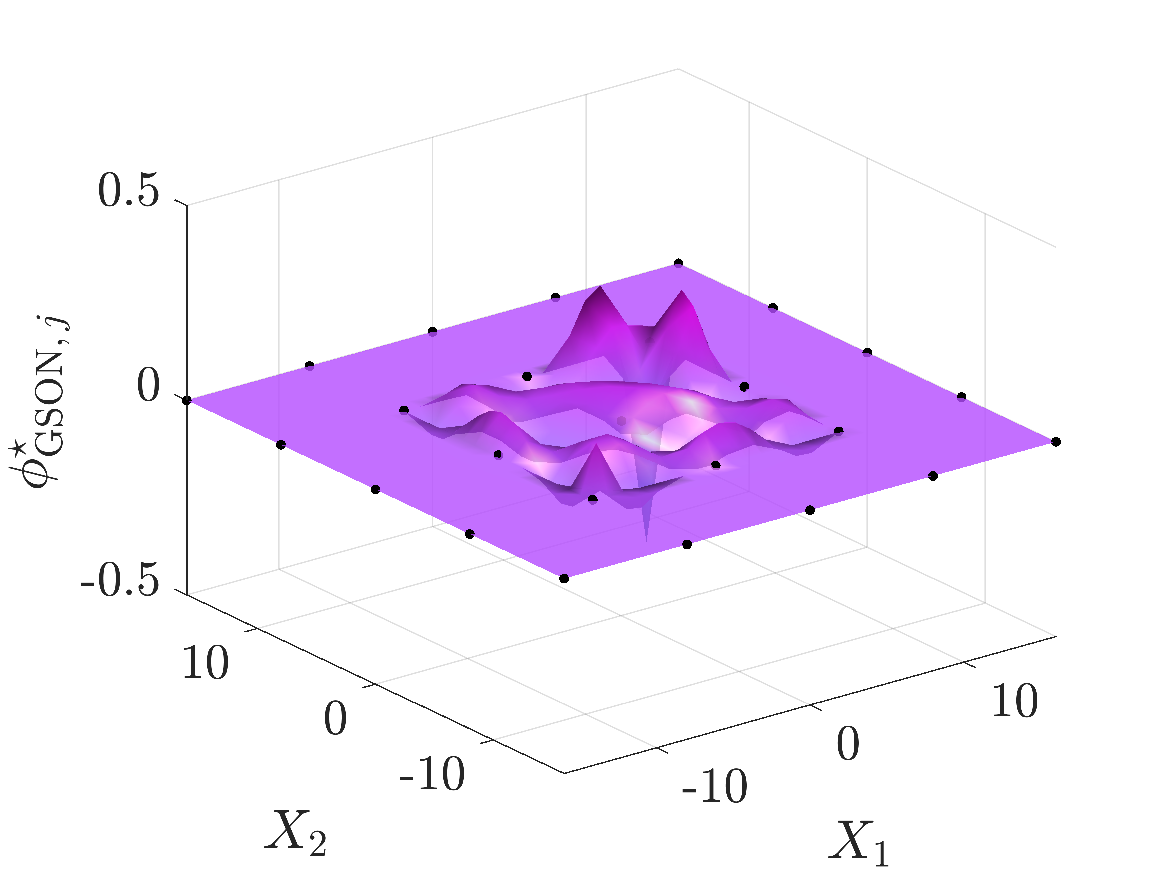}}
   \caption{Enriched LME functions using the shifted Heaviside sign function and orthonormalization with the Gram--Schmidt procedure with $\gamma = 0.9$ for  (a) the center repatom,  (b) a boundary repatom, and (c) a corner repatom of the two dimensional domain. The repatoms are indicated by black dots}
   \label{fig:Schematic_IntFuncHeavGSON_gamma09}
\end{figure}

Multiplying the shifted Heaviside function $\chi_{\mathrm{sign}}(\mathbf{r}_0^{\alpha}) - \chi_{\mathrm{sign}}(\mathbf{r}_0^{\beta_j})$ by the LME shape function $\phi_{\beta_j}(\mathbf{r}_0^{\alpha})$ of the enriched repatom $\beta_j$
yields the enriched shape function $\phi^{\star}_j(\mathbf{r}_{0}^{\alpha})$ as defined in Eq.~(\ref{eq:interpoaltion}).
Figure~\ref{fig:Schematic_IntFuncHeav_gamma09} visualizes $\phi^{\star}_j$ for the three repatoms used in
Figure~\ref{fig:Schematic_IntFunc_gamma09}.
Linear or nearly linear dependence among the enrichment functions leads to ill-conditioning of the discrete system. To stabilize the formulation, we orthonormalize all $\phi_j^{\star}$ using the modified Gram--Schmidt procedure~\cite{Beese2018}. The procedure starts by calculating the orthogonal Gram--Schmidt basis functions for the enriched repatom $j$
\begin{equation}
    \phi_{\mathrm{GS}, \, j}^{\star} 
    = 
    \phi_j^{\star} - \sum_{i=1}^{j-1} 
    \left( 
    \phi_j^{\star \, T} \cdot \phi_{\mathrm{GSON}, \, i}^{\star} 
    \right) 
    \phi_{\mathrm{GSON}, \, i}^{\star}
    \label{eq:GS_basis_func}
\end{equation}
by subtracting a sum term from $\phi_{j}^{\star}$ of all enriched repatoms in the range from $1$ to $j-1$. In the final step, the orthogonalized enriched interpolation functions $\phi_{\mathrm{GS}, \, j}^{\star}$, obtained via the Gram--Schmidt procedure, are normalized as
\begin{equation}
    \phi_{\mathrm{GSON}, \, j}^{\star}
    = 
    \dfrac{\phi_{\mathrm{GS}, \, j}^{\star}}{||\phi_{\mathrm{GS}, \, j}^{\star}||},
    \label{eq:GSON_basis_func}
\end{equation}
resulting in orthonormalized Gram--Schmidt--enriched basis functions $\phi_{\mathrm{GSON}, \, j}^{\star}$ for repatom $j$, as illustrated in Figure \ref{fig:Schematic_IntFuncHeavGSON_gamma09}. 
While orthonormalization removes the direct graphical interpretation of the functions, it preserves weak Kronecker-delta properties, improves the conditioning of the system matrix, and effectively captures discontinuities.

Beyond weak discontinuities with an area, we also model line-type weak discontinuities as fibers and represent them by an enrichment function in the interpolation (see Figure \ref{fig:Schematic_NumExamples}. Following Radtke et al.~\cite{Radtke2010,Radtke2010b,Radtke2011} and Pike and Oskay~\cite{Pike2015}, we use the Heaviside step function $\chi_{\mathrm{step}}$, derived from the signed-distance function $\psi$. $\chi_{\mathrm{step}}$ takes the value $+0.5$ on the fiber (interface) and $0$ elsewhere. The derivation of the fiber-enrichment function is provided in Werner et al.~\cite{Werner2024e7415} and is omitted here for brevity. The orthonormalization procedure described above applies equally to $\chi_{\mathrm{step}}$.

%%%%%%%%%%%%%%%%%%%%%%%%%%%%%%%%%%%%%%%%%%%%%%%%%%%%%%%%%%%%%%%%%%%%%%%%%%%%%
%% Summation rules in QC
%%%%%%%%%%%%%%%%%%%%%%%%%%%%%%%%%%%%%%%%%%%%%%%%%%%%%%%%%%%%%%%%%%%%%%%%%%%%%
\subsection{Summation rules in QC}

Summation (sampling) rules constitute the second pillar of the QC method. They approximate the total potential energy by evaluating only a selected subset of sites or interactions with appropriate weights, thereby reducing the cost of lattice calculations. A common choice is the cluster summation rule~\cite{KNAP20011899,EIDEL200987}, where each repatom is associated with a small cluster of sampling atoms so that the rule is exact for uniform deformation and accurate for low-order deformation gradients. Yang et~al.~\cite{YANG201326} combined cluster sampling with quadrature-point sampling to balance accuracy and cost across non-uniform meshes. The central summation rule was introduced by Beex et al.~\cite{BEEX2014242}, in which repatoms are used together with an additional sampling atom in the center of each interpolation triangle. Roko\v{s} et~al. \cite{ROKOS2017769} used the central summation rule for enriched interpolations to account for strong discontinuities.  Within each subcell of an element created by the crack path, the standard central rule is applied by selecting sampling atoms at the corners and center of a subcell. 

Amelang and co-workers~\cite{AMELANG2015378,AMELANG2015166} provided optimal rules for non-uniform meshes in one, two, and three dimensions. They analyzed error sources and provided practical guidance for selecting cluster sizes and sampling patterns to achieve target accuracy at minimal cost. Kraschewski et~al.~\cite{Kraschewski20241432} used the optimal summation rule for quadratic interpolation in beam lattices, where the energy of an element is represented by sampling unit cells on the vertices, edge centers, and quadrature points of an element. Beex et~al.~\cite{BEEX201552} sampled the interactions directly to approximate the potential energy in an X-braced lattice. One link in each of the four directions (a unit cell) is selected to represent all links within a triangle. Werner et~al.~\cite{Werner2024e7415} extended this interaction-based summation rule to enriched, nonconforming meshes with weak discontinuities by expanding the selected set to include interactions adjacent to the interface, which are sampled individually. The meshless QC framework by Venturini and Kochmann~\cite{Kochmann2014034007} employed a summation rule based on tetrahedral barycenters and Voronoi-weighted volumes, providing a potential foundation for future LME-compatible summation schemes.

A dedicated summation rule for LME interpolation with enrichment has not yet been developed and is beyond the scope of this study. As a pragmatic starting point, we propose building on the summation rule of Venturini and Kochmann~\cite{Kochmann2014034007}, which combines energy evaluations at tetrahedral barycenters with Voronoi-based volume weights, and extending it by discretely sampling all interactions adjacent to interfaces~\cite{Werner2024e7415}.

%%%%%%%%%%%%%%%%%%%%%%%%%%%%%%%%%%%%%%%%%%%%%%%%%%%%%%%%%%%%%%%%%%%%%%%%%%%%%
%% Optimization of the LME Locality Parameter
%%%%%%%%%%%%%%%%%%%%%%%%%%%%%%%%%%%%%%%%%%%%%%%%%%%%%%%%%%%%%%%%%%%%%%%%%%%%%
\subsection{Optimization of the LME locality parameter}
The objective is to determine a set of locality parameters $\boldsymbol{\beta}^{\mathrm{LME}}$ of the LME interpolation functions that minimizes the total potential elastic energy $\Pi(\boldsymbol{\beta}^{\mathrm{LME}})$ of the discrete lattice system and to find kinematically the best LME interpolation. To this end, the following minimization problem is formulated
\begin{equation}
    \boldsymbol{\beta}^{\mathrm{LME}*} = \mathrm{arg \, min}_{\boldsymbol{\beta}^{\mathrm{LME}} \in [\beta^{\mathrm{LME}}_{\min}, \beta_{\max}^{\mathrm{LME}}]} \Pi\left(\boldsymbol{\beta}^{\mathrm{LME}}\right),
    \label{eq:Minimization_Pi_beta}
\end{equation}
with a lower $\beta_{\min}^{\mathrm{LME}} = \gamma^{\mathrm{LME}}_{\min}/ h^2$ and an upper limit $\beta_{\max}^{\mathrm{LME}} = \gamma^{\mathrm{LME}}_{\max}/ h^2$, both defined in terms of the dimensionless locality parameter $\gamma^{\mathrm{LME}}$ and the repatom spacing $h$. The range of $\gamma^{\mathrm{LME}}$ is constrained by practical considerations: if $\gamma^{\mathrm{LME}}$ becomes too small, the stiffness matrix is fully populated and loses its sparsity, making the system of equations increasingly expensive to solve. In this study, the lower limit $\gamma^{\mathrm{LME}}_{\min}$ ranges between $0.3$ and $1.0$. With the upper limit $\gamma^{\mathrm{LME}}_{\max} = 4.0$ it is assumed that the LME interpolation is nearly linear (recall Figure \ref{fig:Schematic_IntFunc_gamma40}).

The minimization problem in Eq. (\ref{eq:Minimization_Pi_beta}) is solved using a quasi-Newton solver using MATLAB's fmincon function with the interior-point and LBFGS algorithm, which requires the gradient of the total elastic energy of the system with respect to $\boldsymbol{\beta}^{\mathrm{LME}}$. Limiting ourselves to  displacement boundary conditions only, the gradient is defined as
\begin{equation}
    \dfrac{\partial \Pi}{\partial \beta^{\mathrm{LME}}_{\beta}} 
    = 
    \left( \dfrac{\partial \boldsymbol{\Phi}}{\partial \beta^{\mathrm{LME}}_{\beta}} \mathbf{r}_{\mathrm{rep}} \right)^T
    \mathbf{G}_{\mathrm{int}},
    \label{eq:dPi_dbeta}
\end{equation}
where the derivative of the interpolation matrix $\boldsymbol{\Phi}$ with respect to $\boldsymbol{\beta}^{\mathrm{LME}} = [\beta^{\mathrm{LME}}_1, \beta^{\mathrm{LME}}_2, \dots, \beta^{\mathrm{LME}}_{n_{\mathrm{rep}}}] \in \mathbb{R}^{N_{\mathrm{rep}}}$ can be obtained by differentiating Eqs. (\ref{eq:LME_shapefunction}), (\ref{eq:GS_basis_func}), and (\ref{eq:GSON_basis_func}), depending on the interpolation scheme. This derivative is multiplied by the location of the repatoms $\mathbf{r}_{\mathrm{rep}}$ and the internal forces at the repatoms $\mathbf{G}_{\mathrm{int}}$.

Using the Gram--Schmidt process described in Section~\ref{subsec:enrichment} has the advantage that the
orthonormalized columns are expressed as a sequence of linear projections and normalizations, and the corresponding
derivatives follow directly from the chain rule. However, Gram--Schmidt is order-dependent and the first enriched
function is only normalized, whereas subsequent functions are successively projected and thus modified more strongly.
Consequently, the ordering influences $\partial \boldsymbol{\Phi}/\partial \beta^{\mathrm{LME}}_{\beta}$ and
$\partial \Pi/\partial \beta^{\mathrm{LME}}_{\beta}$ in Eq.~(\ref{eq:dPi_dbeta}), and ultimately the optimizer’s
result $\boldsymbol{\beta}^{\mathrm{LME}\,*}$ in Eq. (\ref{eq:Minimization_Pi_beta}). Appendix \ref{appendix:app2} describes the computation of the potential energy gradient with respect to $\beta^{\mathrm{LME}}_{\beta}$ and the derivatives of both the LME shape functions and the enriched interpolation functions in more detail.

%%%%%%%%%%%%%%%%%%%%%%%%%%%%%%%%%%%%%%%%%%%%%%%%%%%%%%%%%%%%%%%%%%
%% Numerical Examples
%%%%%%%%%%%%%%%%%%%%%%%%%%%%%%%%%%%%%%%%%%%%%%%%%%%%%%%%%%%%%%%%%%
\section{Numerical Examples}\label{sec3}

This section compares four strategies for selecting the LME locality parameter $\boldsymbol{\gamma}^{\mathrm{LME}}$: the literature-based \emph{baseline}, a globally optimized \emph{uniform} value, a fully optimized \emph{nonuniform} field, and the proposed \emph{pattern-based} approximation. The comparison is carried out for three benchmark problems with weak discontinuities and focuses on the interpolation error and the resulting displacement accuracy relative to the fully resolved reference solutions. Particular attention is paid to whether the optimized locality-parameter fields exhibit systematic spatial structure near material interfaces and whether this structure can be exploited to define simple, non-optimized rules. In this way, the numerical examples are used not only to assess accuracy, but also to identify a practical alternative to computationally expensive per-repatom optimization.

\begin{figure}[tbp]
   \centering
   \subfloat[Circular Inclusion]{\includegraphics[width=0.32\textwidth]{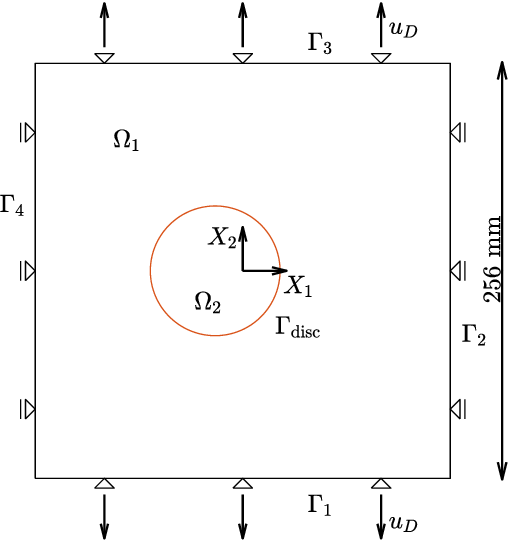}}
   \hfill
   \subfloat[Square Inclusion]{\includegraphics[width=0.32\textwidth]{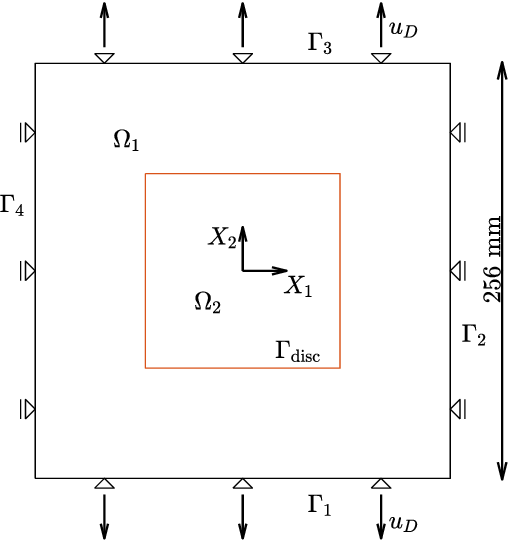}}
   \hfill
   \subfloat[Fiber]{\includegraphics[width=0.32\textwidth]{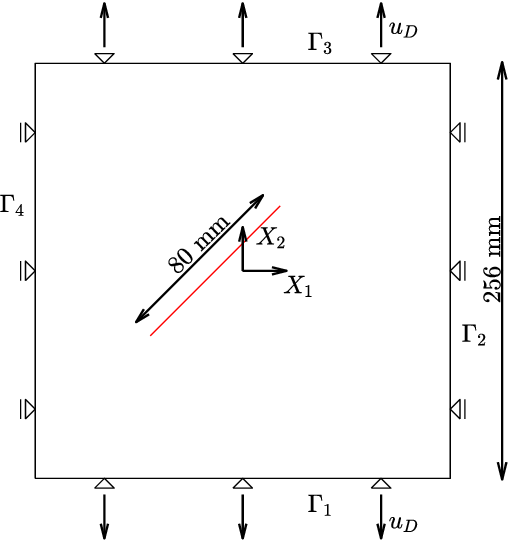}}
   \hfill
   \caption{Three numerical examples (a) circular inclusion, (b) square inclusion, and (c) fiber along with the applied boundary conditions, including the prescribed displacement $u_D$ on $\Gamma_1$ and $\Gamma_3$.}
   \label{fig:Schematic_NumExamples}
\end{figure}

To study the influence of the LME locality parameter with enriched interpolation and its precision, we consider three numerical examples (Figure~\ref{fig:Schematic_NumExamples}). 
Each example uses a square domain of $\SI{256}{mm}\times\SI{256}{mm}$ with lattice spacing $d=\SI{1}{mm}$. The opposite boundaries $\Gamma_1$ and $\Gamma_3$ are subjected to a prescribed displacement $u_{\mathrm D}$. All four boundaries are fixed in the $X_1$ direction, while atoms and repatoms on $\Gamma_2$ and $\Gamma_4$ are unconstrained in the $X_2$ direction.

The circular inclusion in Fig.~\ref{fig:Schematic_NumExamples}~a has a radius $\SI{40}{mm}$ and a center point $\mathbf{P}_{\mathrm{cent}}=(\SI{-17}{mm},\,\SI{0}{mm})$. The square inclusion in Fig.~\ref{fig:Schematic_NumExamples}~b has an edge length of $\SI{60}{mm}$ and is centered in the domain. The fiber in Fig.~\ref{fig:Schematic_NumExamples}~c is $\SI{80}{mm}$ long and centered at $\mathbf{P}_{\mathrm{cent}}=(\SI{-17}{mm},\,\SI{0}{mm})$. The Young's modulus of the truss elements within the inclusions is $10$ times that of the surrounding matrix; for the fiber it is $100$ times. For reference, in heterogeneous materials such as concrete, aggregates typically exhibit a Young’s modulus about four times that of the cement matrix, \cite{Xue20231134} whereas reinforcing fibers may be roughly twelve times stiffer.\cite{GUO2021113072}

These three benchmarks were selected to evaluate the interpolation under distinct interface geometries and singularity types. (i) The \emph{circular inclusion} introduces a curved interface whose normal varies with respect to the underlying X-braced lattice. The interface has a polygonal representation on the lattice following horizontal and vertical links that produce many local corners. (ii) The \emph{square inclusion} has piecewise straight interfaces aligned with the underlying lattice and exhibits corners and moderate stress singularities, enabling evaluation of the interpolation without curvature. (iii) The \emph{fiber} provides a single straight interface, with the two fiber tips generating the strongest localized fields (tip singularities) as a result of the high stiffness contrast.

The optimum support size of the LME interpolation is problem-dependent \cite{Rosolen2010868} and therefore an entire range of the dimensionless locality parameter has been applied in the literature. Arroyo and Ortiz \cite{Arroyo20062167} suggested $\gamma^{\mathrm{LME}} = 1.8$ as optimal for non-linear elastic problems. Rosolen et al. \cite{Rosolen2010868} used $\gamma^{\mathrm{LME}} = 1.6$ as a reference and compared it with optimized $\boldsymbol{\gamma}^{\mathrm{LME}}$ fields for hyperelastic materials. Kochmann and Venturini \cite{Kochmann2014034007} assumed $\gamma^{\mathrm{LME}} = 1.0$ as optimal for LME interpolation for QC.

Given the diversity of “optimal” locality settings reported for different applications, our goal is to determine suitable $\gamma^{\mathrm{LME}}$ values for heterogeneous lattices with Heaviside-enriched interpolation and nonconforming parameterization of the domain. We extend the optimization strategy of Rosolen et al.~\cite{Rosolen2010868} to the enriched basis and apply it to three inclusion benchmarks. First, we optimize a single uniform $\gamma^{\mathrm{LME}}$, which is the same for all interpolation functions and repatoms. Second, we optimize a spatially varying field $\boldsymbol{\gamma}^{\mathrm{LME}}$ and analyze the resulting patterns, particularly near the interfaces. We identify simple default rules that deliver accuracy comparable to that of the optimized field without additional optimization.

For this purpose, we use seven different interpolation schemes with LME or linear interpolation, and with and without Heaviside enrichment: 
\begin{enumerate}
  \item \emph{xQC LME Baseline--H:}
        LME with $\gamma^{\mathrm{LME}}=1.8$ (Arroyo \& Ortiz~\cite{Arroyo20062167}) and Heaviside enrichment. 
        This represents the current state of the art for a uniform LME locality parameter.
  \item \emph{xQC LME Optimized--Uniform--H:}
        LME + Heaviside enrichment; a single locality parameter $\gamma^{\mathrm{LME}}$ is optimized and applied uniformly to all interpolation functions.
  \item \emph{xQC LME Optimized--Nonuniform--H:}
        LME + Heaviside enrichment; a spatially varying field $\boldsymbol{\gamma}^{\mathrm{LME}}$ is obtained by optimizing the locality parameter per repatom  and is nonuniform in space.
  \item \emph{xQC LME Optimized--Nonuniform (no H):}
        LME without Heaviside enrichment; as \emph{xQC LME Optimized--Nonuniform--H}, optimize $\boldsymbol{\gamma}^{\mathrm{LME}}$ per repatom and nonuniform in space, but without enrichment.
  \item \emph{xQC LME Pattern--Based--H:}
        Spatially varying $\boldsymbol{\gamma}^{\mathrm{LME}}$ assigned by simple rules extracted from \emph{xQC LME Optimized--Nonuniform--H} and \emph{xQC LME Optimized--Nonuniform} without running an optimizer.
  \item \emph{sQC Linear:}
        Standard QC with linear interpolation and a conforming mesh leading to a fully resolved interface.
  \item \emph{xQC Linear--H:}
        Extended QC with linear interpolation and Heaviside enrichment and therefore a nonconforming mesh.
\end{enumerate}

\begin{figure}[tbp]
   \centering
   \subfloat[xQC LME]{\includegraphics[width=0.32\textwidth]{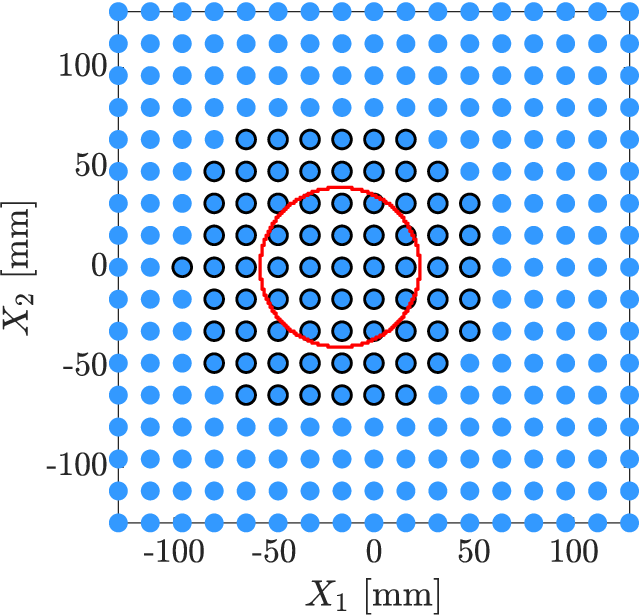}}
   \hfill
   \subfloat[sQC Linear]{\includegraphics[width=0.32\textwidth]{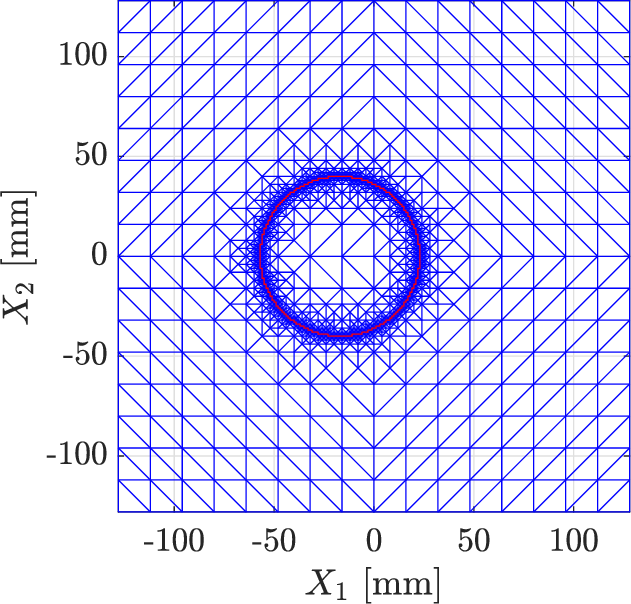}}
   \hfill
   \subfloat[xQC Linear--H]{\includegraphics[width=0.32\textwidth]{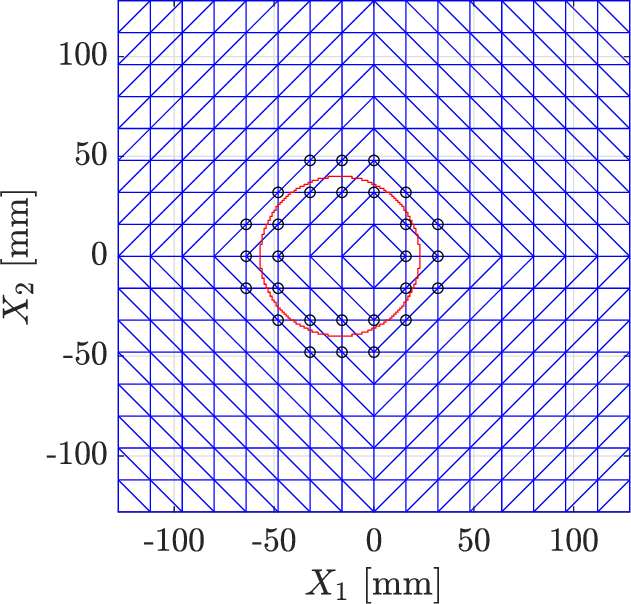}}
   \hfill
   \caption{Comparison of interpolation strategies: (a) xQC–LME, (b) sQC–linear, and (c) xQC–linear.}
   \label{fig:Discretization_NumExamples}
\end{figure}

For LME interpolation schemes with Heaviside enrichment (1-3 and 5 from the above list) the repatoms within a distance of $2.5h$ from the interface are enriched for circular and square inclusions, and within a distance of $0.7h$ for the fiber example. Figure \ref{fig:Discretization_NumExamples}~a shows representatively the circular inclusion example with enriched repatoms (black circles). 
In addition to the LME interpolation schemes with Heaviside enrichment, the examples are also solved with the standard QC method (\emph{sQC Linear}) using linear interpolation and a fully resolved interface (Figure \ref{fig:Discretization_NumExamples}~b). The extended QC method using linear interpolation combined with Heaviside enrichment is denoted by \emph{xQC Linear--H}, in which repatoms adjacent to the interface are enriched (Figure~\ref{fig:Discretization_NumExamples}~c). The interpolation schemes \emph{sQC Linear} and \emph{xQC Linear--H} are included for comparison purposes to provide context for the results, and are discussed in more detail in Werner et al\cite{Werner2024e7415}.

To assess the accuracy of individual interpolation schemes, the full solution always serves as a reference. The number of DOFs for all the systems considered is summarized in Table~\ref{tab:DOF_CircIncl} for the circular inclusion example. The other two examples differ in the exact numbers of DOFs but remain comparable and hence not shown. The fully resolved solution involves approximately \SI{133000}{} unknowns, while the QC approaches achieve one order of magnitude reduction.
\begin{table}[tbp]
\centering
\caption{Degrees of freedom (DOFs) for the QC interpolation schemes at various repatom spacings $h$ for the circular inclusion example (full solution \SI{133000}{} unknowns).}
\begin{tabular}{lrrrrr}
\hline
Method & 32 mm & 16 mm & 8 mm & 4 mm & 2 mm \\
\hline
xQC LME       & 250  & 726  & 2494 & 9102   & -- \\
sQC Linear    & 3150   & 3442   & 4818   & 10658   & 34674 \\
xQC Linear--H & 162   & 578   & 2178   & 8450   & 33282 \\
\hline
\end{tabular}
\label{tab:DOF_CircIncl}
\end{table}

Accuracy of all QC schemes is assessed in terms of the relative displacement error quantified as
\begin{align}
\varepsilon^{\mathbf{u}} = \frac{\normtwo{\mathbf{u}^{\mathrm{QC}} - \mathbf{u}^{\mathrm{FS}} } }{\normtwo{ \mathbf{u}^{\mathrm{FS}} }},
\label{eq:Error_Disp}
\end{align}
where $\mathbf{u}^{\mathrm{QC}}$ denotes the global displacements obtained from the QC analysis, while $\mathbf{u}^{\mathrm{FS}}$ corresponds to the reference full solution, in both cases $\mathbf{u} = \mathbf{r} - \mathbf{r}_0$. We report only the relative displacement error, since other quantities such as the relative elastic energy or stress error exhibit similar trends and lead to similar conclusions. In addition to the global displacement error, the local per-atom error in the displacement is defined as
\begin{equation}
  \varepsilon^{\mathbf{u}}_{\alpha}
  = \bigl|\, \normtwo{ \mathbf{u}^{\mathrm{QC}}_{\alpha} }
          - \normtwo{ \mathbf{u}^{\mathrm{FS}}_{\alpha} } \,\bigr|
          \label{eq:AbsError}
\end{equation}
and will be reported.

%%%%%%%%%%%%%%%%%%%%%%%%%%%%%%%%%%%%%%%%%%%%%%%%%%%%%%%%%%%%%%%%%%
%% Uniformly Optimized Locality Parameter
%%%%%%%%%%%%%%%%%%%%%%%%%%%%%%%%%%%%%%%%%%%%%%%%%%%%%%%%%%%%%%%%%%
\subsection{Baseline and uniform optimized locality parameter}
\label{sec:UniformlyOptimized}

To evaluate the performance of the baseline value of $\gamma^{\mathrm{LME}} = 1.8$, as suggested by Arroyo and Ortiz \cite{Arroyo20062167}, we first optimize this parameter and then compare the accuracy of the two interpolation schemes, \emph{xQC LME Baseline–H} and \emph{xQC LME Optimized–Uniform–H}. The optimized $\gamma^{\mathrm{LME}}$ values for the circular and square inclusion examples range between 1.26 and 1.91 (Table \ref{tab:gamma_uniform}), which are close to the literature value of $\gamma^{\mathrm{LME}} = 1.8$. 
\begin{table}[tbp]
\centering
\caption{Optimized uniform dimensionless locality parameter $\gamma^{\mathrm{LME}}$ for the interpolation scheme \emph{xQC LME Optimized--Uniform--H} of the three numerical examples with four different repatom spacings.}
\begin{tabular}{lcccc}
\hline
& \multicolumn{4}{c}{$\gamma^{\mathrm{LME}}$ at $h$ [mm]} \\
\cline{2-5}
Numerical example & 32 & 16 & 8 & 4 \\
\hline
Circular inclusion & 1.26 & 1.58 & 1.34 & 1.64 \\
Square inclusion   & 1.28 & 0.80 & 1.27 & 1.91 \\
Fiber              & 0.80 & 0.80 & 0.80 & 0.80 \\
\hline
\end{tabular}
\label{tab:gamma_uniform}
\end{table}
An exception to this range occurs for the square inclusion with a repatom spacing of $h = \SI{16}{mm}$, for which the optimizer consistently identified the minimum at the lower bound, regardless of the initial value of $\gamma^{\mathrm{LME}}$.

The fiber example yields a constant $\gamma^{\mathrm{LME}} = 0.8$ across all repatom spacings. These differences in $\gamma^{\mathrm{LME}}$ are likely related to the different enrichment strategies in the respective domains: the circular and square inclusion examples include enriched repatoms up to a distance of $2.5h$ from the interface, whereas the fiber example uses a shorter enrichment distance of $0.7h$. The two enrichment distances used here were determined by iterative refinement to avoid high condition number of the system of equations and numerical instability. A more systematic approach to selecting the enrichment distance could be developed, but this is considered beyond the scope of the present contribution.

\begin{figure}[tbp]
  \captionsetup{font=normalsize}
  \centering
  \includegraphics[width=\textwidth]{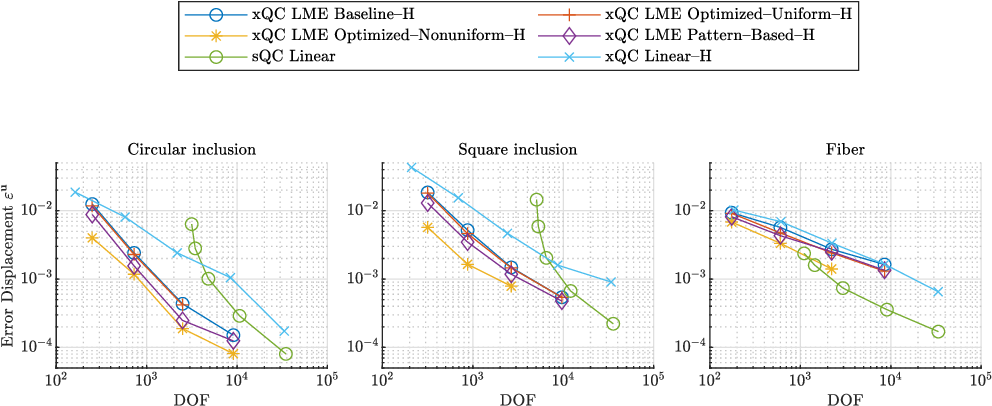}
  \caption{Relative error in displacement for the three numerical examples: circular inclusion (left), square inclusion (center), and fiber (right). Six interpolation schemes are applied to each example.}
  \label{fig:ErDispDOF_all_three}
\end{figure}

The relative error in displacement (Eq.~\ref{eq:Error_Disp}) results in six curves shown as double logarithmic plots in Figure~\ref{fig:ErDispDOF_all_three}. Here, we focus on the interpolation schemes \emph{xQC LME Optimized–Uniform–H} and \emph{xQC LME Baseline–H}; the others are discussed in the following subsections.
The interpolation scheme \emph{xQC LME Optimized--Uniform--H} leads only to slightly more accurate displacements compared to \emph{xQC LME Baseline--H}. However, both schemes yield \SIrange{17}{63}{\percent} of the \emph{xQC Linear--H} error, with the greatest improvements at finer repatom spacings. In the fiber example, \emph{xQC LME Optimized--Uniform--H} offers only a slight improvement over \emph{xQC Linear--H}, likely because the fiber–matrix Young’s modulus contrast is 100:1, far exceeding the 10:1 contrast of the circular and square inclusions, and fewer enriched repatoms.

% \emph{xQC LME Baseline--H}
% \emph{xQC LME Optimized--Uniform--H}
% \emph{xQC LME Optimized--Nonuniform--H}
% \emph{xQC LME Optimized--Nonuniform (no H)}
% \emph{xQC LME Pattern--Based--H}
% \emph{sQC Linear}
% \emph{xQC Linear--H}

%%%%%%%%%%%%%%%%%%%%%%%%%%%%%%%%%%%%%%%%%%%%%%%%%%%%%%%%%%%%%%%%%%
%% Nonuniform Optimized Locality Parameter
%%%%%%%%%%%%%%%%%%%%%%%%%%%%%%%%%%%%%%%%%%%%%%%%%%%%%%%%%%%%%%%%%%
\subsection{Nonuniform optimized locality parameter}
\label{sec:NonuniformOptimized}

Nonuniform (per repatom) optimization of $\boldsymbol{\gamma}^{\mathrm{LME}}$ is computationally expensive and not sustainable for large discrete lattice models. We therefore optimize the LME locality parameter to identify patterns of $\boldsymbol{\gamma}^{\mathrm{LME}}$ for weak discontinuities along curved interfaces with many corners (circular inclusion), along straight interfaces (square inclusion), and in regions with high stress concentration (fiber). For this purpose we consider two interpolation schemes, \emph{xQC LME Optimized--Nonuniform--H} and \emph{xQC LME Optimized--Nonuniform (no H)}.

\begin{figure}[tbp]
   \centering
   \subfloat[Circular Inclusion]{\includegraphics[width=0.32\textwidth]{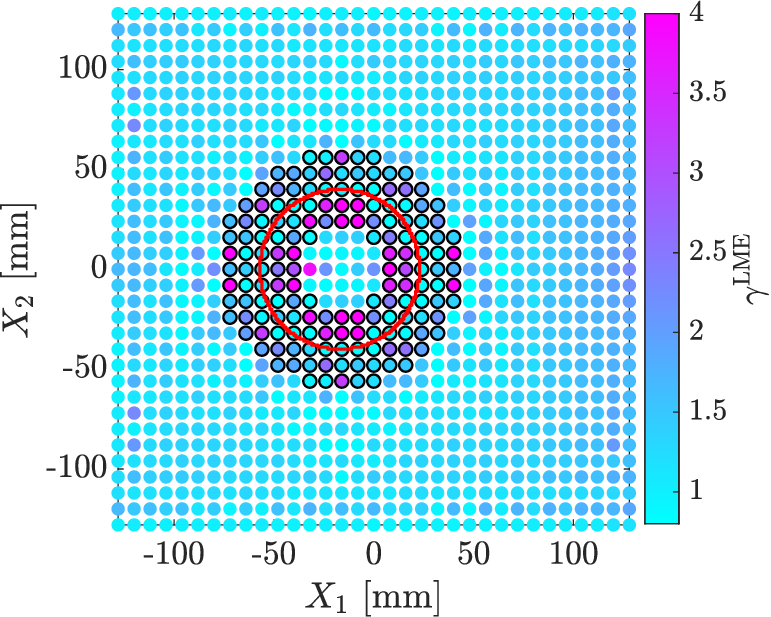}}
   \hfill
   \subfloat[Square Inclusion]{\includegraphics[width=0.32\textwidth]{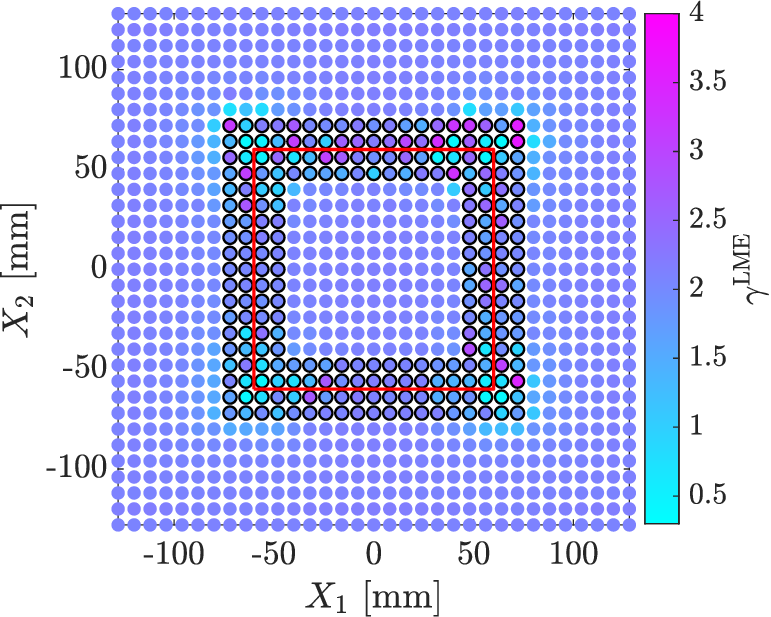}}
   \hfill
   \subfloat[Fiber]{\includegraphics[width=0.32\textwidth]{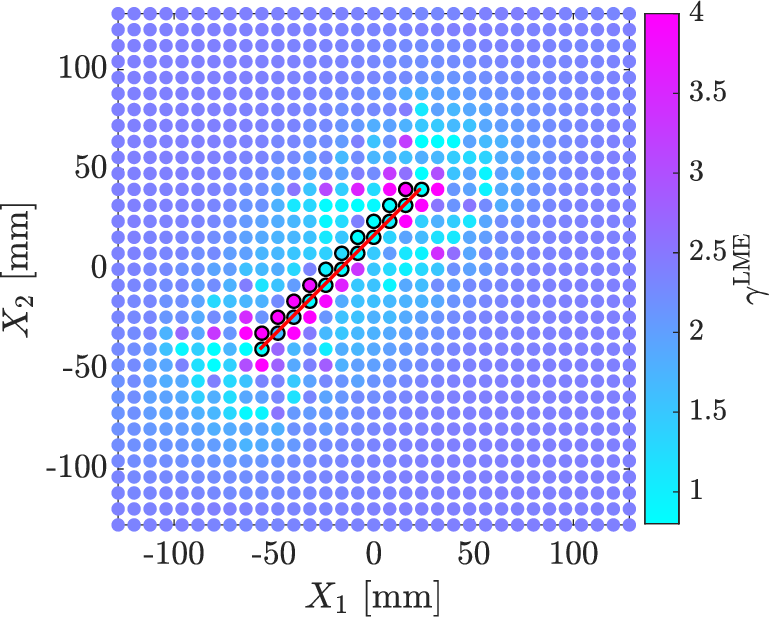}}
   \caption{Distribution of $\boldsymbol{\gamma}^{\mathrm{LME}}$ using the \emph{xQC LME Optimized--Nonuniform--H} interpolation scheme with a repatom spacing $h=\SI{8}{mm}$ for the (a) circular inclusion, (b) square inclusion, and (c) fiber.}
   \label{fig:GammaDistr}
\end{figure}

In Figure \ref{fig:GammaDistr}, the spatial distribution of the optimized dimensionless locality parameter $\boldsymbol{\gamma}^{\mathrm{LME}}$ is shown for the three examples with a repatom spacing $h=\SI{8}{mm}$ using the interpolation scheme \emph{xQC LME Optimized--Nonuniform--H}. The enriched repatoms are indicated by black circles. The optimization of $\boldsymbol{\gamma}^{\mathrm{LME}}$ is performed with a lower bound of $0.3$ for the square inclusion, a lower bound of $0.8$ for the circular inclusion and the fiber, and an upper bound of $4.0$ for all three of them.
In all three examples, the far-field values of $\boldsymbol{\gamma}^{\mathrm{LME}}$ increase gradually (circular inclusion in Figure \ref{fig:GammaDistr}~a) or remain nearly constant between 2.0 and 2.3 (square inclusion and fiber in Figure \ref{fig:GammaDistr}~b~and~c). In contrast, the repatoms located near the interface, as well as the enriched repatoms, exhibit considerable scatter. Furthermore, no symmetric patterns of the dimensionless locality parameter are observed, which is attributed to the order dependence of the Gram--Schmidt procedure.

\begin{figure}[tbp]
   \centering
   \subfloat[$\gamma^{\mathrm{LME}}$ Distribution]{\includegraphics[width=0.32\textwidth]{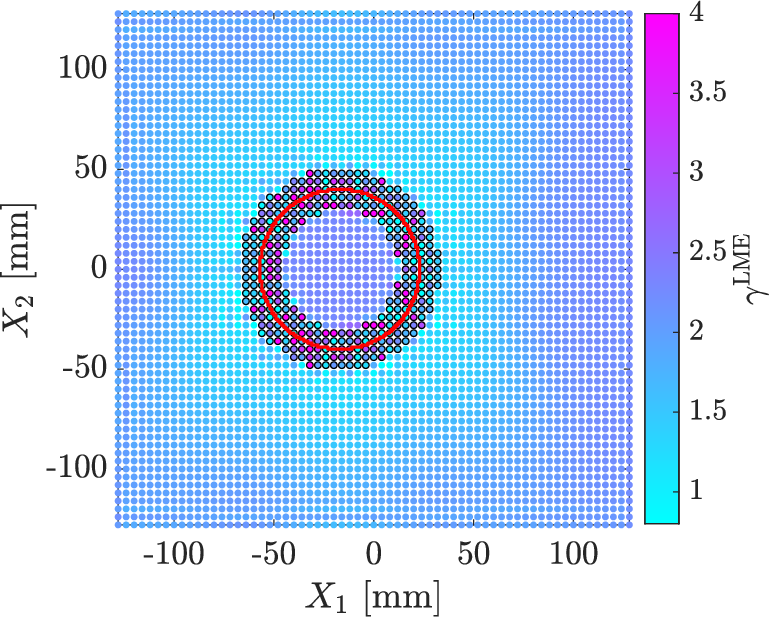}}
   \hfill
   \subfloat[$\gamma^{\mathrm{LME}}$ vs. Signed Distance]{\includegraphics[width=0.32\textwidth]{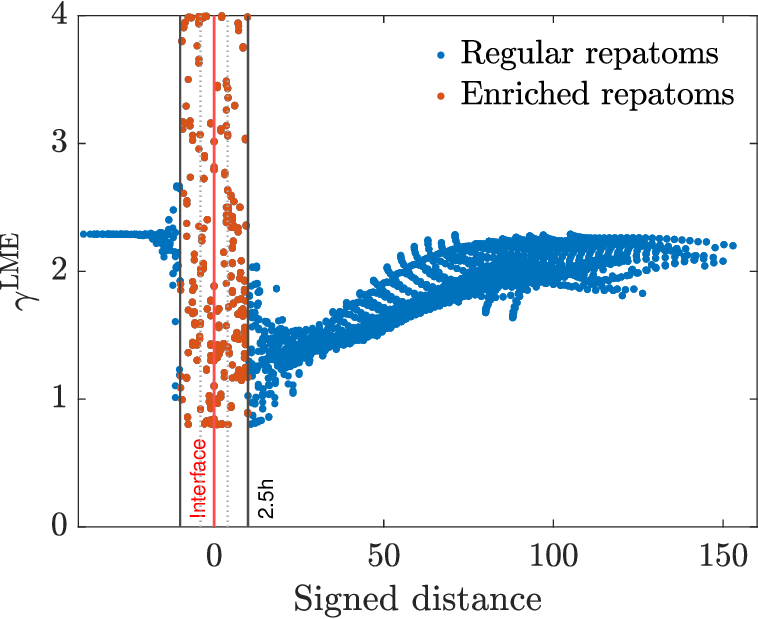}}
   \hfill
   \subfloat[Median with IQR band]{\includegraphics[width=0.32\textwidth]{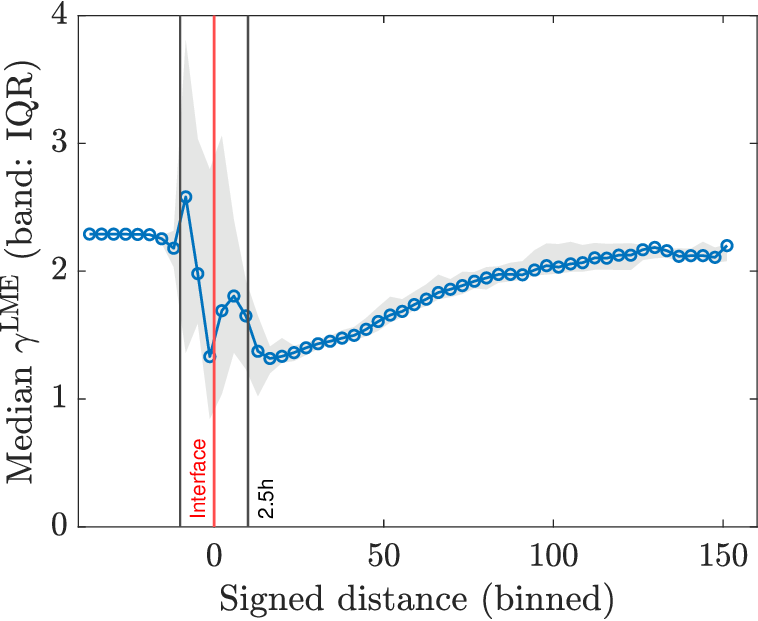}}
   \caption{(a) Distribution of $\boldsymbol{\gamma}^{\mathrm{LME}}$ for the circular inclusion, $h=\SI{4}{mm}$, using the interpolation scheme \emph{xQC LME Optimized--Nonuniform--H}; (b) $\boldsymbol{\gamma}^{\mathrm{LME}}$ at each repatom as a function of the signed distance (from the interface); (c) bin-wise interquartile range (IQR) of $\boldsymbol{\gamma}^{\mathrm{LME}}$ with the corresponding median.}
   \label{fig:GammaDistr_CircIncl_h4}
\end{figure}

Figure~\ref{fig:GammaDistr_CircIncl_h4}a shows the resulting $\boldsymbol{\gamma}^{\mathrm{LME}}$ field for the circular inclusion analyzed with the \emph{xQC LME Optimized–Nonuniform–H} interpolation scheme and a repatom spacing of $h = \SI{4}{mm}$.
In Figure~\ref{fig:GammaDistr_CircIncl_h4}~b, the $\boldsymbol{\gamma}^{\mathrm{LME}}$ values for regular and enriched repatoms are plotted as a function of the signed distance from the interface. The interface itself is indicated by a red line at zero distance, with negative values corresponding to repatoms located inside the inclusion. The enrichment zone with a width of $2.5h$ is highlighted by two vertical black lines. Figure~\ref{fig:GammaDistr_CircIncl_h4}~c presents the bin-wise interquartile range (IQR) of $\boldsymbol{\gamma}^{\mathrm{LME}}$ together with the corresponding medians. Repatoms are grouped into 54 bins according to their signed distance, and for each bin, the 25th and 75th percentiles (Q1 and Q3) define the shaded IQR interval [Q1, Q3], while the bin-wise median is shown in blue.

Within the enriched zone, $\boldsymbol{\gamma}^{\mathrm{LME}}$ varies between approximately 0.8 and 4.0, confirming the pronounced scatter without any apparent dependence on the signed distance. Inside the inclusion, $\boldsymbol{\gamma}^{\mathrm{LME}}$ remains nearly constant at around 2.3, whereas in the far field, the values increase gradually. For regular repatoms with positive signed distances, $\boldsymbol{\gamma}^{\mathrm{LME}}$ rises smoothly from about 1.4 to 2.2. The median value in Figure~\ref{fig:GammaDistr_CircIncl_h4}~c indicated that the bin closest to the interface has a relatively low $\boldsymbol{\gamma}^{\mathrm{LME}}$, and increases in the two adjacent bins on either side within the enriched repatom zone. Nevertheless, the IQR remains wide, reflecting the pronounced variability of $\boldsymbol{\gamma}^{\mathrm{LME}}$ in the enriched repatom zone.

The observed nonuniform distribution of $\boldsymbol{\gamma}^{\mathrm{LME}}$ can be attributed to the Heaviside enrichment and, in particular, to the order-dependent orthonormalization of the enriched basis functions in the Gram--Schmidt process (see Equation~\ref{eq:GS_basis_func}). By construction, the orthonormalization introduces asymmetries in the enriched shape functions, leading to differences in the interpolation of repatoms that are otherwise symmetry-equivalent in axisymmetric or point-symmetric configurations. The subsequent optimization compensates for these variations by adjusting the locality parameters, which manifests as a nonuniform field of $\boldsymbol{\gamma}^{\mathrm{LME}}$ (Figure~\ref{fig:GammaDistr_CircIncl_h4}~a).

\begin{figure}[tbp]
   \centering
   \subfloat[$\gamma^{\mathrm{LME}}$ Distribution]{\includegraphics[width=0.32\textwidth]{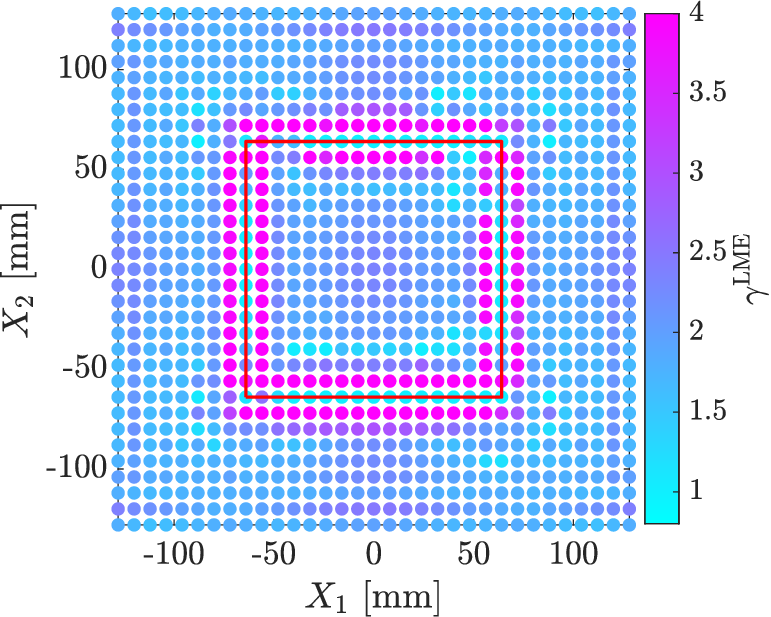}}
   \hfill
   \subfloat[$\gamma^{\mathrm{LME}}$ vs. Signed Distance]{\includegraphics[width=0.32\textwidth]{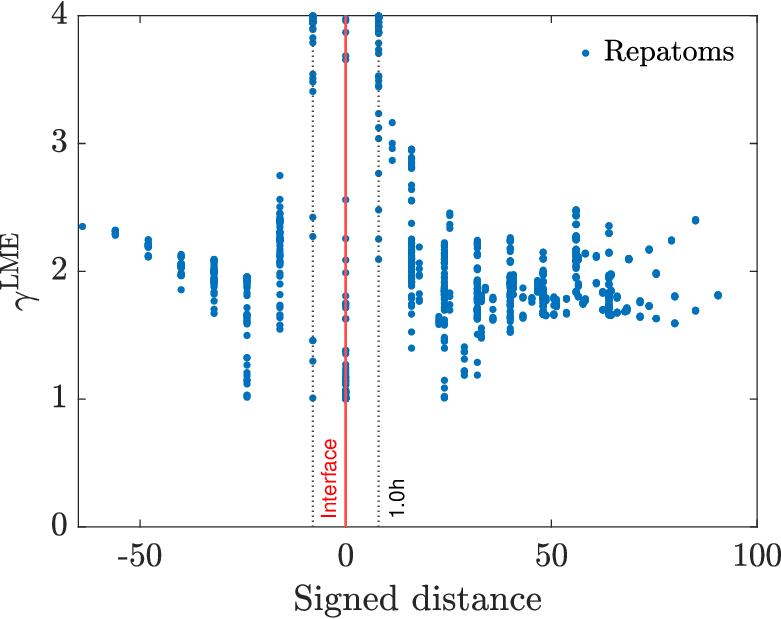}}
   \hfill
   \subfloat[Median with IQR band]{\includegraphics[width=0.32\textwidth]{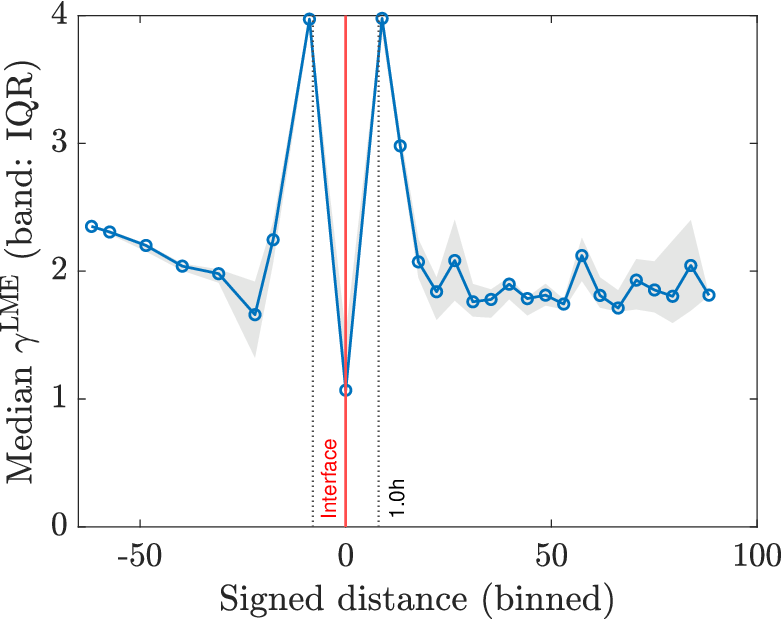}}
   \caption{(a) Distribution of $\boldsymbol{\gamma}^{\mathrm{LME}}$ for the modified square inclusion and $h=\SI{8}{mm}$, with the interface aligned with the repatoms and using the interpolation scheme \emph{xQC LME Optimized–Nonuniform (no H)}; (b) $\boldsymbol{\gamma}^{\mathrm{LME}}$ at each repatom as a function of the signed distance (from the interface); (c) bin-wise interquartile range (IQR) of $\boldsymbol{\gamma}^{\mathrm{LME}}$ with the corresponding median.}
   \label{fig:GammaDistr_SquareIncl_woEnr}
\end{figure}

To better understand the influence of Heaviside enrichment and the order-dependent orthonormalization in the Gram--Schmidt process, a modified version of the square inclusion example is analyzed using the \emph{xQC LME Optimized–Nonuniform (no H)} interpolation scheme. The inclusion is slightly larger, extending from $\SI{-64}{mm}$ to $\SI{64}{mm}$ in the $X_1$ and $X_2$ directions, and aligns with the repatoms (Figure \ref{fig:GammaDistr_SquareIncl_woEnr}~a). 
A repatom spacing of $h = \SI{8}{mm}$ is used, and the optimization is carried out with lower and upper bounds of 1.0 and 4.0 for $\boldsymbol{\gamma}^{\mathrm{LME}}$. In Figures \ref{fig:GammaDistr_SquareIncl_woEnr}~a~and~b, a distinct pattern of the dimensionless locality parameter $\boldsymbol{\gamma}^{\mathrm{LME}}$ can already be observed, which becomes even more pronounced in Figure \ref{fig:GammaDistr_SquareIncl_woEnr}~c, where the bin-wise IQR and median values are shown. The median on the interface $\boldsymbol{\gamma}^{\mathrm{LME}} = 1.07$ is close to the lower bound while one repatom spacing away from both sides of the interface, the median $\boldsymbol{\gamma}^{\mathrm{LME}} = 3.97$ is close to the upper bound.

\begin{figure}[tbp]
   \centering
   \subfloat[$\gamma^{\mathrm{LME}}$ Distribution]{\includegraphics[width=0.32\textwidth]{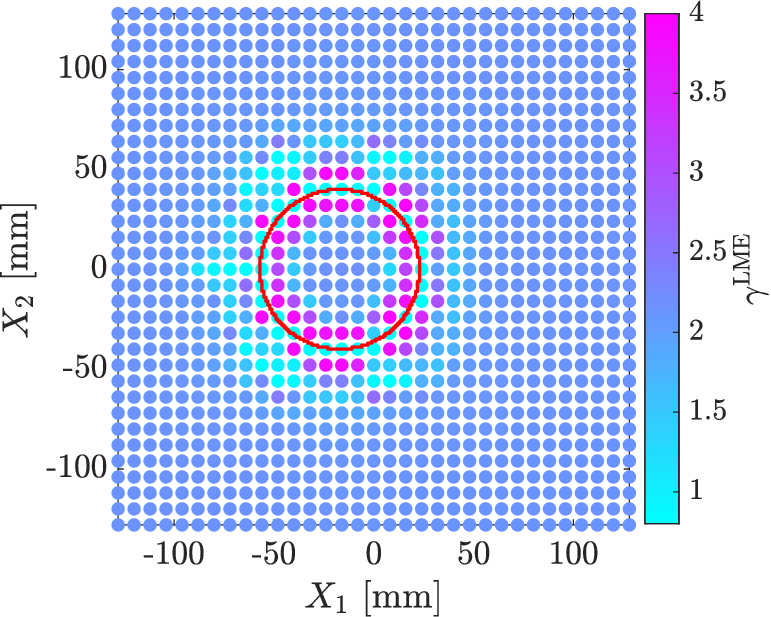}}
   \hfill
   \subfloat[$\gamma^{\mathrm{LME}}$ vs. Signed Distance]{\includegraphics[width=0.32\textwidth]{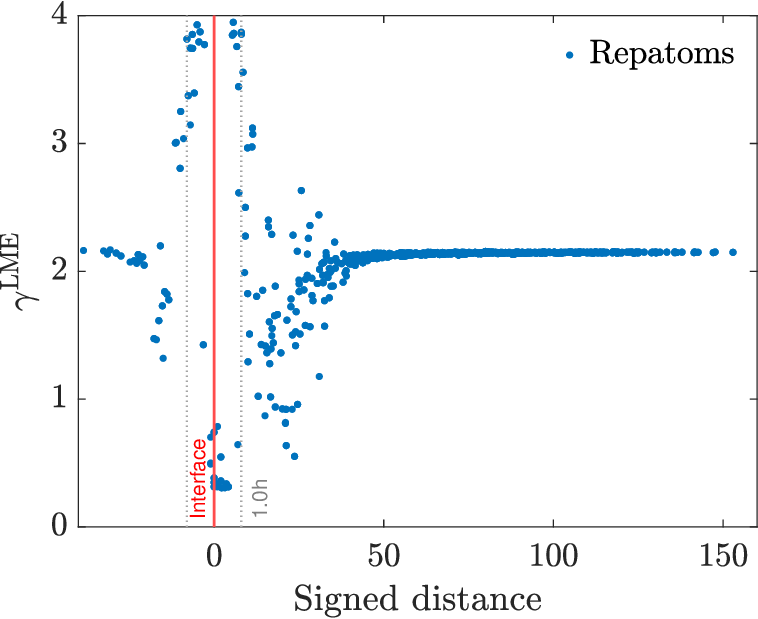}}
   \hfill
   \subfloat[Median with IQR band]{\includegraphics[width=0.32\textwidth]{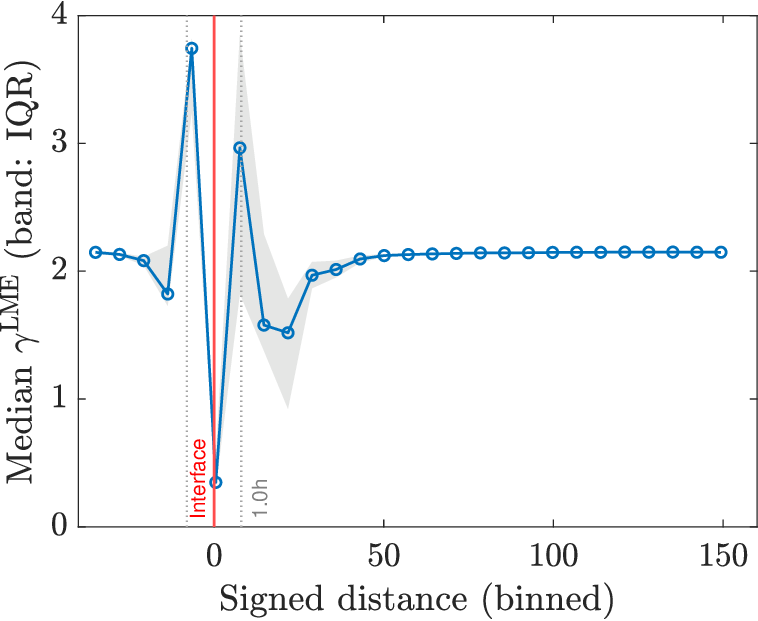}}
   \caption{(a) Distribution of $\boldsymbol{\gamma}^{\mathrm{LME}}$ for the circular inclusion, a repatom spacing of $h=\SI{8}{mm}$, and using the interpolation scheme \emph{xQC LME Optimized–Nonuniform (no H)}; (b) $\boldsymbol{\gamma}^{\mathrm{LME}}$ at each repatom as a function of the signed distance (from the interface); (c) bin-wise interquartile range (IQR) of $\boldsymbol{\gamma}^{\mathrm{LME}}$ with the corresponding median.}
   \label{fig:GammaDistr_CircIncl_woEnr}
\end{figure}

In addition to the modified square inclusion, the circular inclusion with a repatom spacing of $h=\SI{8}{mm}$ is analyzed using the \emph{xQC LME Optimized–Nonuniform (no H)} interpolation scheme (Figure~\ref{fig:GammaDistr_CircIncl_woEnr}). 
Unlike the modified square case, the interface cannot be aligned with the repatom grid. The values of $\boldsymbol{\gamma}^{\mathrm{LME}}$ near the interface tend toward the lower bound, here set to $\boldsymbol{\gamma}^{\mathrm{LME}}_{\min} = 0.3$ (Figure~\ref{fig:GammaDistr_CircIncl_woEnr}~b–c). The median value of the dimensionless locality parameter in the bin closest to the interface is approximately 0.35. In the far field, $\boldsymbol{\gamma}^{\mathrm{LME}}$ approaches an almost constant value of approximately 2.15. Overall, the characteristic \emph{low–high–intermediate} trend of $\boldsymbol{\gamma}^{\mathrm{LME}}$ with increasing distance from the interface, observed for the modified square inclusion, can also be identified for the circular inclusion analyzed using the \emph{xQC LME Optimized–Nonuniform (no H)} interpolation scheme. Consequently, the \emph{xQC LME Optimized–Nonuniform (no H)} interpolation scheme exhibits clear and reproducible patterns (Figures \ref{fig:GammaDistr_SquareIncl_woEnr} and \ref{fig:GammaDistr_CircIncl_woEnr}), distinct from the behavior observed for the \emph{xQC LME Optimized--Nonuniform--H} scheme (Figure~\ref{fig:GammaDistr_CircIncl_h4}).

The interpolation scheme \emph{xQC LME Optimized–Nonuniform (no H)} is excluded from the accuracy assessment, as it is used solely to analyze the distribution patterns of $\boldsymbol{\gamma}^{\mathrm{LME}}$. Due to the absence of Heaviside enrichment, this scheme naturally exhibits lower accuracy in reproducing displacements. In contrast, \emph{xQC LME Optimized–Nonuniform–H} results in the smallest relative displacement errors for the circular and square inclusions (Figure~\ref{fig:ErDispDOF_all_three}). Its errors are \SI{32}{\percent}–\SI{54}{\percent} of those of \emph{xQC LME Optimized–Uniform–H} and \SI{8}{\percent}–\SI{22}{\percent} of those of \emph{xQC Linear–H}. For the fiber example, \emph{xQC LME Optimized–Nonuniform–H} leads to approximately one-half of the \emph{xQC Linear–H} error and about two-thirds of the \emph{xQC LME Optimized–Uniform–H} error. Compared to the inclusion cases, the fiber example appears less sensitive to the adaptation of the support-size near the interface, which explains the more modest gain of the LME-based variants.

% \emph{xQC LME Baseline--H}
% \emph{xQC LME Optimized--Uniform--H}
% \emph{xQC LME Optimized--Nonuniform--H}
% \emph{xQC LME Optimized--Nonuniform (no H)}
% \emph{xQC LME Pattern--Based--H}
% \emph{sQC Linear}
% \emph{xQC Linear--H}

%%%%%%%%%%%%%%%%%%%%%%%%%%%%%%%%%%%%%%%%%%%%%%%%%%%%%%%%%%%%%%%%%%
%% Pattern-Based Locality Parameter
%%%%%%%%%%%%%%%%%%%%%%%%%%%%%%%%%%%%%%%%%%%%%%%%%%%%%%%%%%%%%%%%%%
\subsection{Pattern-based locality parameter}
\label{sec:PatternBased}

Based on the analysis carried out in Section~\ref{sec:NonuniformOptimized}, distinct patterns can be identified in the nonuniform optimized $\boldsymbol{\gamma}^{\mathrm{LME}}$ fields (Figures~\ref{fig:GammaDistr_SquareIncl_woEnr} and~\ref{fig:GammaDistr_CircIncl_woEnr}). In this section, we quantify these patterns and propose an explicit analytical formulation for $\boldsymbol{\gamma}^{\mathrm{LME}}$ in the vicinity of the interface. We further demonstrate that, when combined with enrichment, the proposed formulation can achieve an accuracy comparable to or exceeding that of $\boldsymbol{\gamma}^{\mathrm{LME}}$ distributions reported in the literature, while avoiding computationally expensive optimization.

Without enrichment, $\boldsymbol{\gamma}^{\mathrm{LME}}$ exhibits a pronounced dependence on the signed distance, following a characteristic \emph{low–high–intermediate} profile: small values near the interface, a rapid increase across a transition region, and a plateau or gradual rise farther away. In the far field, $\boldsymbol{\gamma}^{\mathrm{LME}}$ approaches a nearly constant value of about $\boldsymbol{\gamma}^{\mathrm{LME}}_{\mathrm{FF}} \approx 2$, while near the interface it tends toward the lower bound imposed by the minimization problem.

\begin{figure}[tbp]
   \centering
   \subfloat[$h=32$]{\includegraphics[width=0.33\textwidth]{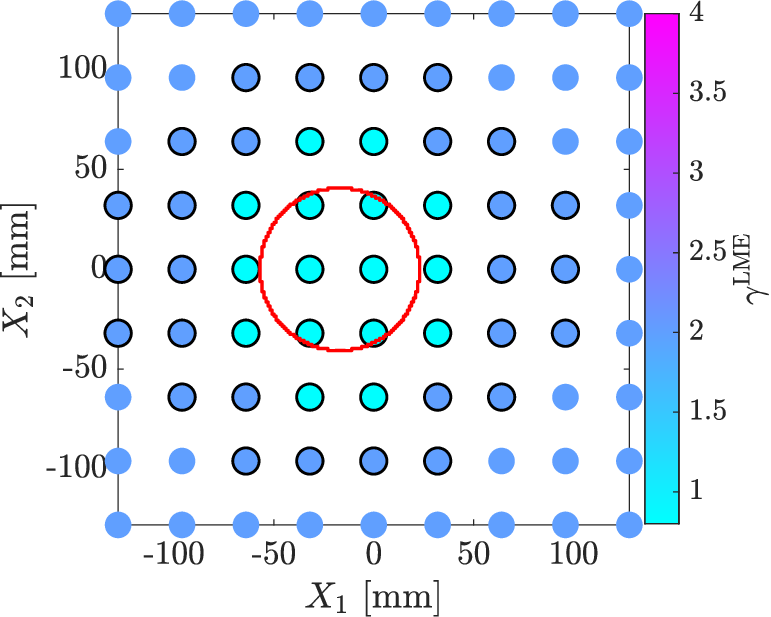}}
   \hfill
   \subfloat[$h=8$]{\includegraphics[width=0.33\textwidth]{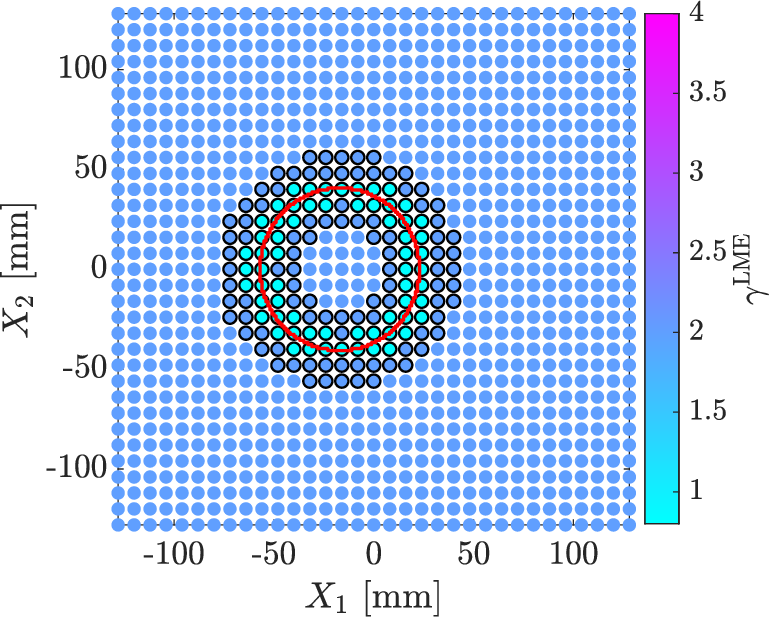}}
   \hfill
   \subfloat[$h=4$]{\includegraphics[width=0.33\textwidth]{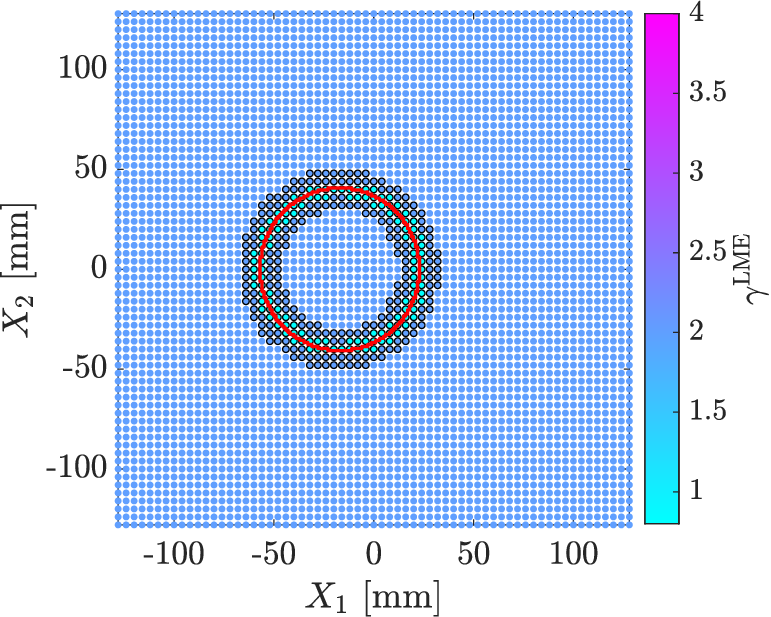}}
   \caption{$\boldsymbol{\gamma}^{\mathrm{LME}}$ distribution defined by the \emph{xQC LME Pattern–Based–H} interpolation scheme for the circular inclusion example with three repatom spacings, $h = [32, 8, 4]$. Enriched repatoms are indicated by black circles. The locality parameter is set to $\boldsymbol{\gamma}^{\mathrm{LME}}_{\mathrm{IF}} = 0.8$ for repatoms within one repatom spacing of the interface and to $\boldsymbol{\gamma}^{\mathrm{LME}}_{\mathrm{FF}} = 2.0$ elsewhere.}
   \label{fig:GammaDistr_CircIncl_PD}
\end{figure}

To capture this behavior in a simple yet representative manner, we define the locality parameter as
\begin{equation}
    \boldsymbol{\gamma}^{\mathrm{LME}} =
    \begin{cases}
        \boldsymbol{\gamma}^{\mathrm{LME}}_{\mathrm{IF}} = 0.8, & \mathrm{for}\; h \leq 1 \\
        \boldsymbol{\gamma}^{\mathrm{LME}}_{\mathrm{FF}} = 2.0, & \mathrm{for}\; h > 1
    \end{cases}
    \label{eq:AnalyticalGamma}
\end{equation}
where $\boldsymbol{\gamma}^{\mathrm{LME}}_{\mathrm{IF}}$ and $\boldsymbol{\gamma}^{\mathrm{LME}}_{\mathrm{FF}}$ represent the locality parameters near the interface and in the far field. For this analytical formulation, the high (peak) part of the characteristic \emph{low–high–intermediate} profile is intentionally omitted (recall Figures \ref{fig:GammaDistr_SquareIncl_woEnr}~c and \ref{fig:GammaDistr_CircIncl_woEnr}~c. Preliminary tests that included this feature led to lower overall accuracy compared to the proposed simplified formulation. Moreover, the \emph{xQC LME Optimized–Nonuniform–H} interpolation scheme does not exhibit the pronounced local maximum (Figure~\ref{fig:GammaDistr_CircIncl_h4}) observed in the \emph{xQC LME Optimized–Nonuniform (no H)} case, supporting the omission of the peak in the analytical definition. The interpolation scheme \emph{xQC LME Pattern–Based–H}, defined by Eq.~(\ref{eq:AnalyticalGamma}), is illustrated in Figure~\ref{fig:GammaDistr_CircIncl_PD} for three repatom spacings in the circular inclusion example.
The case with a repatom spacing of $h = \SI{8}{mm}$ in Figure~\ref{fig:GammaDistr_CircIncl_PD}~b can be directly compared to the corresponding results in Figure~\ref{fig:GammaDistr_CircIncl_woEnr} and Figure~\ref{fig:GammaDistr}~a.

% \emph{xQC LME Baseline--H}
% \emph{xQC LME Optimized--Uniform--H}
% \emph{xQC LME Optimized--Nonuniform--H}
% \emph{xQC LME Optimized--Nonuniform (no H)}
% \emph{xQC LME Pattern--Based--H}
% \emph{sQC Linear}
% \emph{xQC Linear--H}

The interpolation scheme \emph{xQC LME Pattern–Based–H} achieves accuracies comparable to those of \emph{xQC LME Optimized–Uniform–H} in terms of the relative displacement error and exhibits a convergence behavior similar to both \emph{xQC LME Optimized–Uniform–H} and \emph{xQC LME Optimized–Nonuniform–H} (Figure~\ref{fig:ErDispDOF_all_three}). Relative to \emph{xQC Linear–H}, the \emph{xQC LME Pattern–Based–H} scheme yields errors of \SI{10}{\percent}–\SI{50}{\percent} for the circular inclusion, \SI{24}{\percent}–\SI{32}{\percent} for the square inclusion, and \SI{67}{\percent}–\SI{86}{\percent} for the fiber example. Overall, it provides a substantial gain in accuracy while avoiding computationally expensive optimization.

\begin{figure}[tbp]
   \centering
   \subfloat[\emph{xQC LME Optimized–Uniform–H}]{\includegraphics[width=0.33\textwidth]{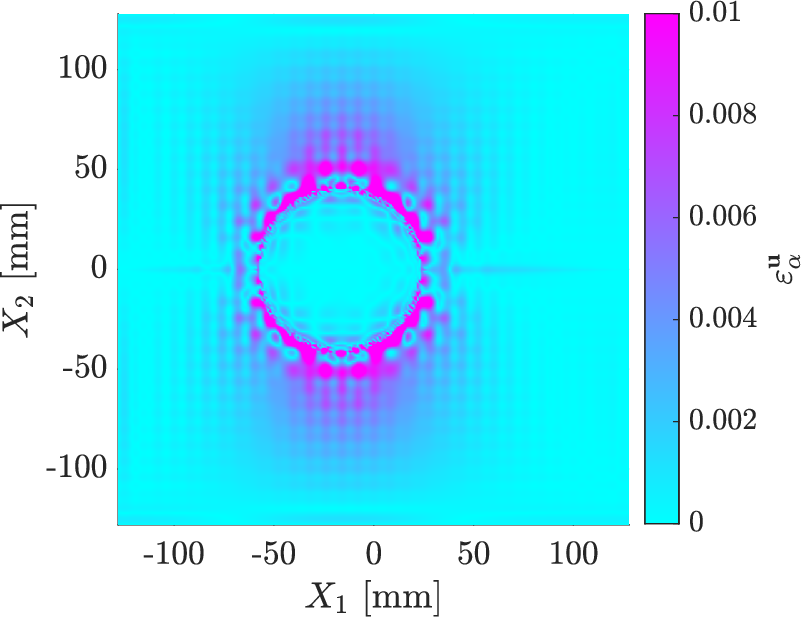}}
   \hfill
%   \subfloat[\emph{xQC LME Optimized–Nonuniform–H}]
   \subfloat[\shortstack{\emph{xQC LME Optimized--}\\\emph{Nonuniform--H}}]{\includegraphics[width=0.33\textwidth]{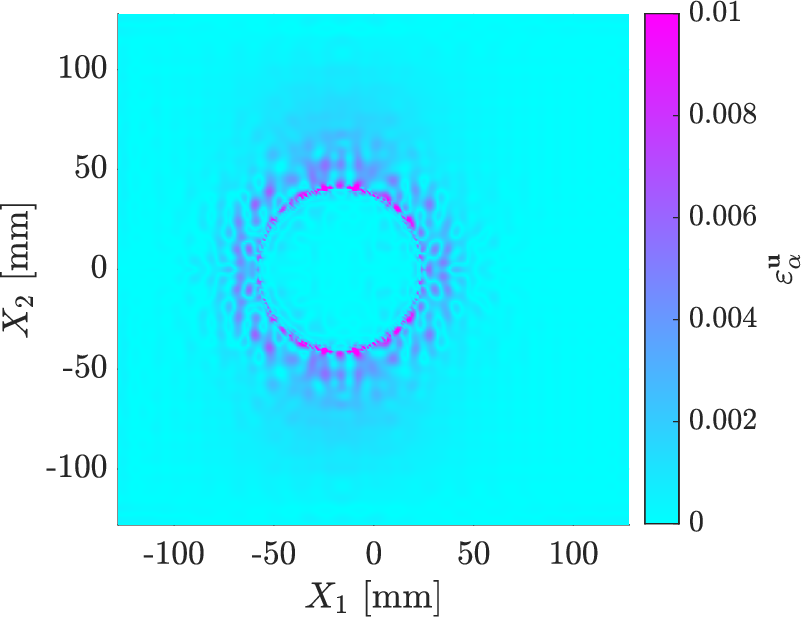}}
   \hfill
   \subfloat[\emph{xQC LME Pattern–Based–H}]{\includegraphics[width=0.33\textwidth]{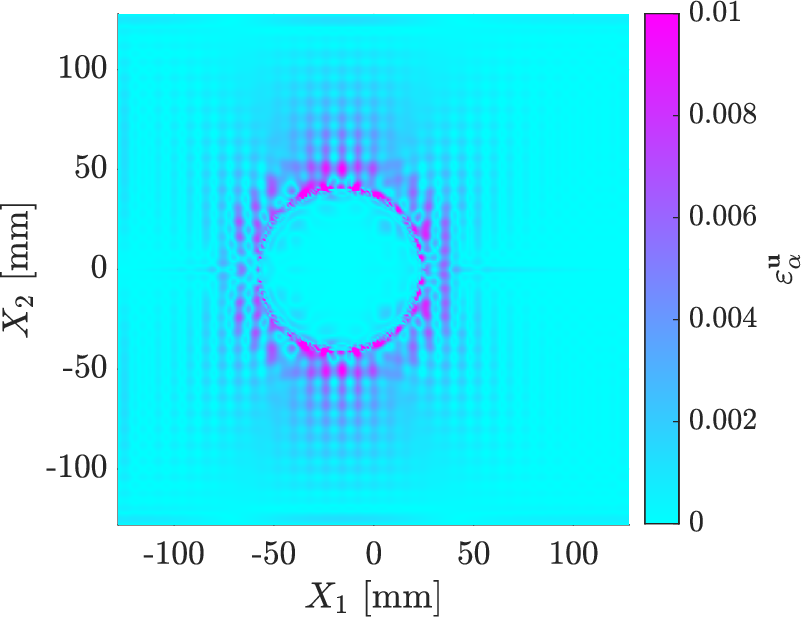}}
   \caption{Local displacement error $\varepsilon^{\mathbf{u}}_{\alpha}$ for the circular inclusion example with a repatom spacing of $h = \SI{8}{mm}$: (a) \emph{xQC LME Optimized–Uniform–H}, (b) \emph{xQC LME Optimized–Nonuniform–H}, and (c) \emph{xQC LME Pattern–Based–H}.}
   \label{fig:AbsErrorDisp_CircIncl}
\end{figure}

For completeness, the displacement error field, computed according to Eq.~(\ref{eq:AbsError}), is shown in Figure~\ref{fig:AbsErrorDisp_CircIncl} for the circular inclusion example with a repatom spacing of $h = \SI{8}{mm}$, comparing the interpolation schemes \emph{xQC LME Optimized–Uniform–H}, \emph{xQC LME Optimized–Nonuniform–H}, and \emph{xQC LME Pattern–Based–H}. The \emph{xQC LME Optimized–Uniform–H} approach (Figure~\ref{fig:AbsErrorDisp_CircIncl}~a) has pronounced errors along the interface and in regions above and below the inclusion, suggesting an influence of the loading direction. In contrast, \emph{xQC LME Optimized–Nonuniform–H} (Figure~\ref{fig:AbsErrorDisp_CircIncl}~b) shows errors concentrated primarily at the interface. With the \emph{xQC LME Pattern–Based–H} scheme (Figure~\ref{fig:AbsErrorDisp_CircIncl}~c), the error regions are clearly reduced compared to \emph{xQC LME Optimized–Uniform–H}.

\begin{figure}[tbp]
   \centering
   \subfloat[\emph{xQC LME Optimized–Uniform–H}]{\includegraphics[width=0.33\textwidth]{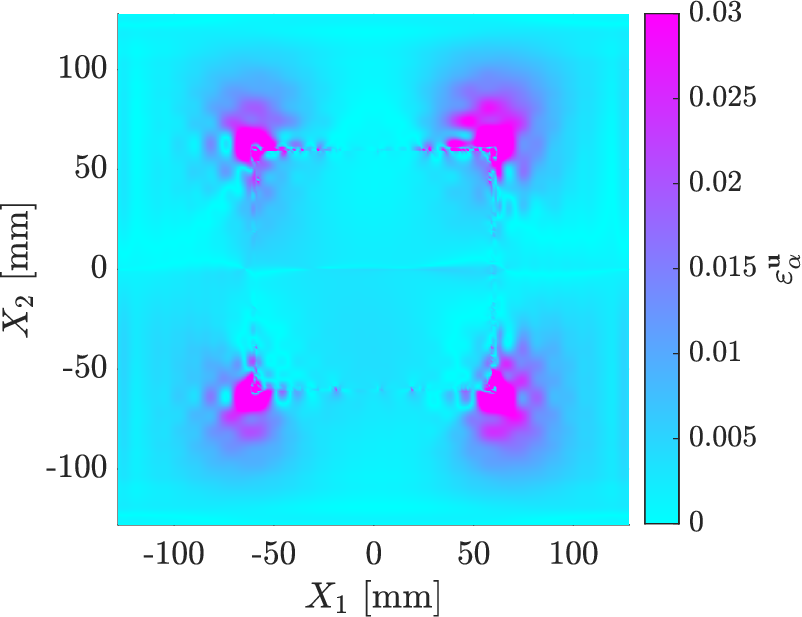}}
   \hfill
   %\subfloat[\emph{xQC LME Optimized–Nonuniform–H}]
   \subfloat[\shortstack{\emph{xQC LME Optimized--}\\\emph{Nonuniform--H}}]{\includegraphics[width=0.33\textwidth]{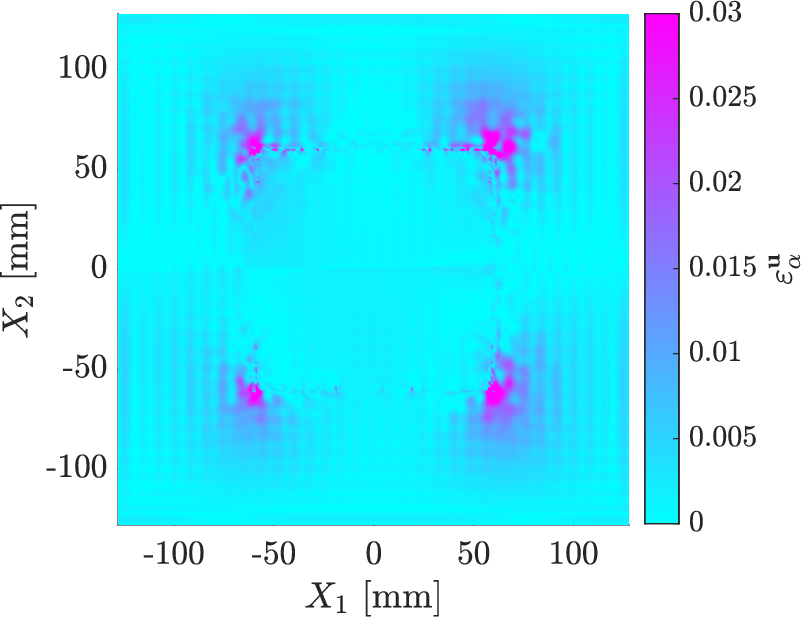}}
   \hfill
   \subfloat[\emph{xQC LME Pattern–Based–H}]{\includegraphics[width=0.33\textwidth]{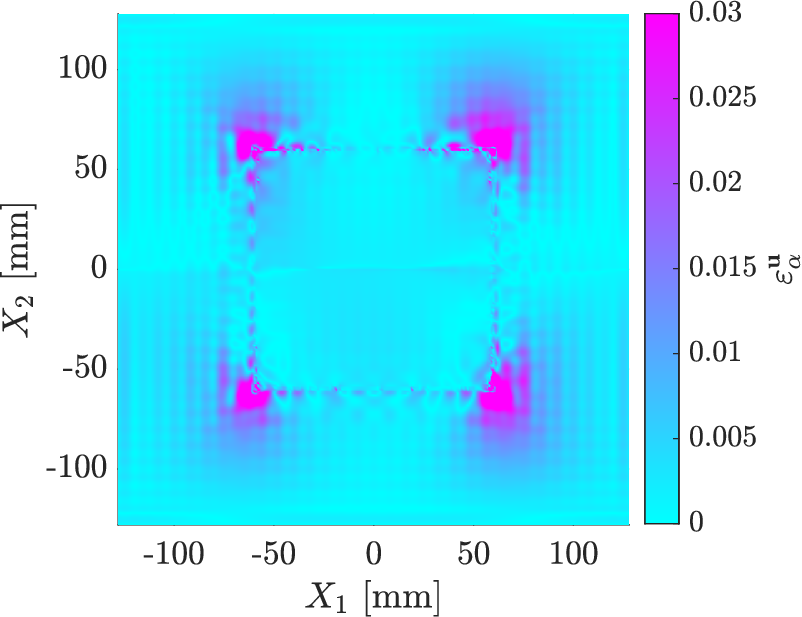}}
   \caption{Local displacement error $\varepsilon^{\mathbf{u}}_{\alpha}$ for the square inclusion example with a repatom spacing of $h = \SI{8}{mm}$: (a) \emph{xQC LME Optimized–Uniform–H}, (b) \emph{xQC LME Optimized–Nonuniform–H}, and (c) \emph{xQC LME Pattern–Based–H}.}
   \label{fig:AbsErrorDisp_SquareIncl}
\end{figure}

Each scheme exhibits displacement errors concentrated at the corners of the square inclusion (Figure~\ref{fig:AbsErrorDisp_SquareIncl}). The \emph{xQC LME Optimized–Nonuniform–H} scheme (b) produces the smallest and most localized error regions, whereas the \emph{xQC LME Optimized–Uniform–H} (a) and \emph{xQC LME Pattern–Based–H} (c) schemes show very similar, more extended error distributions. The observed asymmetry of the error fields is attributed to the order dependence of the orthonormalization in the Gram--Schmidt process.

\begin{figure}[tbp]
   \centering
   \subfloat[\emph{xQC LME Optimized–Uniform–H}]{\includegraphics[width=0.33\textwidth]{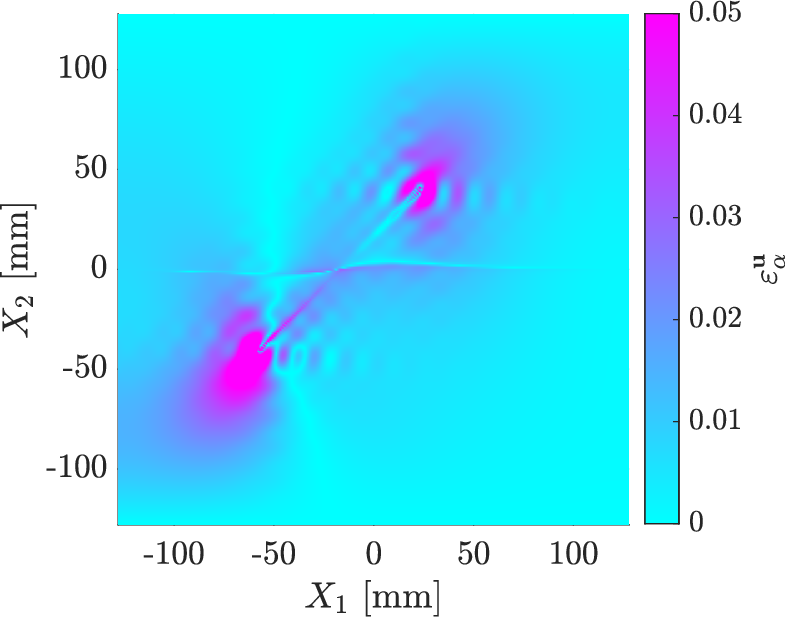}}
   \hfill
  % \subfloat[\emph{xQC LME Optimized–Nonuniform–H}]
   \subfloat[\shortstack{\emph{xQC LME Optimized--}\\\emph{Nonuniform--H}}]{\includegraphics[width=0.33\textwidth]{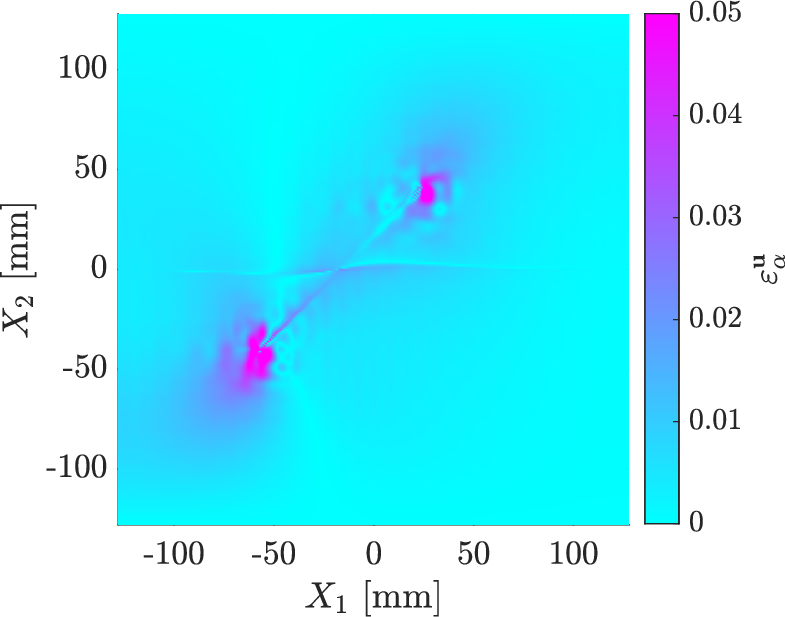}}
   \hfill
   \subfloat[\emph{xQC LME Pattern–Based–H}]{\includegraphics[width=0.33\textwidth]{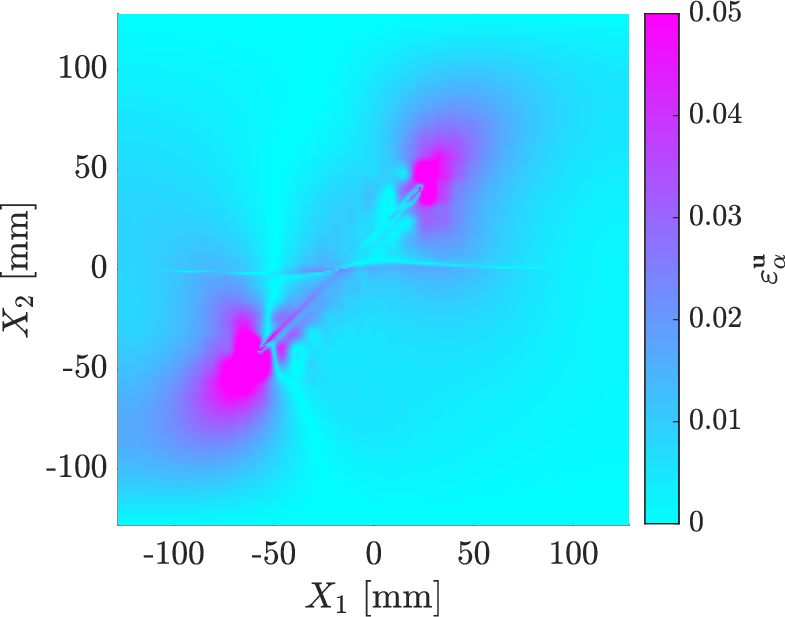}}
   \caption{Local displacement error $\varepsilon^{\mathbf{u}}_{\alpha}$ for the fiber example with a repatom spacing of $h = \SI{8}{mm}$: (a) \emph{xQC LME Optimized–Uniform–H}, (b) \emph{xQC LME Optimized–Nonuniform–H}, and (c) \emph{xQC LME Pattern–Based–H}.}
   \label{fig:AbsErrorDisp_Fiber}
\end{figure}

All three schemes show displacement errors concentrated near the fiber tips (Figure~\ref{fig:AbsErrorDisp_Fiber}), with the \emph{xQC LME Optimized–Nonuniform–H} variant (b) exhibiting the smallest and most localized error regions. In contrast, the \emph{xQC LME Optimized–Uniform–H} (a) and \emph{xQC LME Pattern–Based–H} (c) schemes display broader error zones extending from both fiber ends, consistent with their comparable global relative errors (Figure~\ref{fig:ErDispDOF_all_three}c).

%Across all three examples, the \emph{xQC LME Pattern–Based–H} scheme exhibits local displacement error fields of similar or smaller extent compared to \emph{xQC LME Optimized–Uniform–H}, confirming the accuracy and robustness of the approach. These results support the choice of $\boldsymbol{\gamma}^{\mathrm{LME}}_{\mathrm{IF}} = 0.8$ at the interface and $\boldsymbol{\gamma}^{\mathrm{LME}}_{\mathrm{FF}} = 2.0$ in the far field as a simple and effective default for achieving high precision in enriched LME interpolation without the need for computational optimization.

Across all three benchmark problems, the optimized locality-parameter fields are clearly nonuniform and show their strongest variation in the vicinity of the material interface, whereas the far field is much less sensitive. Although the exact distributions depend on the geometry of the inclusion or fiber, the resulting patterns are sufficiently systematic to motivate a simple rule-based approximation of $\boldsymbol{\gamma}^{\mathrm{LME}}$. In particular, the choice of $\boldsymbol{\gamma}^{\mathrm{LME}}_{\mathrm{IF}} = 0.8$ at the interface and $\boldsymbol{\gamma}^{\mathrm{LME}}_{\mathrm{FF}} = 2.0$ in the far field provides a simple and effective default. These values are inferred from the recurring structure of the optimized nonuniform fields observed across the benchmark cases. The proposed pattern-based field reproduces much of the accuracy gain achieved by the fully optimized nonuniform distributions while avoiding the computational cost of per-repatom optimization. This makes the pattern-based strategy the most relevant outcome of the present study from a practical QC perspective.

\section{Summary and Conclusion}\label{sec4}

In this work, we extend the QuasiContinuum (QC) method by integrating Heaviside enrichment with the Local Maximum Entropy (LME) interpolation to capture weak discontinuities in discrete lattices. We further optimize the LME locality parameter, controlling the basis function support, by minimizing the total elastic potential energy of the system. We consider both global (\emph{uniform}) and per-repatom (\emph{nonuniform}) optimization of the locality parameter and compare the results with the state-of-the-art values of the locality parameter from the literature (\emph{baseline}). The \emph{nonuniform} approach enables the identification of characteristic patterns that can be used to predefine the locality-parameter field. The resulting \emph{pattern-based} approach preserves high accuracy while avoiding computationally expensive optimization. Based on comparative studies performed for multiple single-inclusion and fiber examples, the main conclusions are summarized as follows:

\begin{itemize}
    \item LME interpolation can be readily combined with Heaviside enrichment. A Gram--Schmidt orthonormalization of the enriched basis is suitable to reduce near--linear dependence and improve the conditioning of the resulting system of equations. Using the Heaviside enrichment, the interpolation scheme captures weak discontinuities such as inclusion boundaries and embedded fibers.
    \item Across the three examples, the extended QC variants with LME interpolation (\emph{xQC LME}) yield lower displacement errors than the extended QC with linear basis functions (\emph{xQC Linear}). For the fiber case, the gain is modest, with errors remaining at about \SIrange{42}{89}{\percent} of the \emph{xQC Linear} level. For the circular and square inclusions, the improvements are stronger, with errors dropping to roughly \SIrange{10}{63}{\percent} for the circular inclusion and \SIrange{10}{32}{\percent} for the square inclusion compared to \emph{xQC Linear}.
    \item An explicit \emph{pattern-based} formulation for the nonuniform locality parameter $\boldsymbol{\gamma}^{\mathrm{LME}}$ is proposed as a function of the distance from the interface, improving the accuracy while eliminating the need for costly optimization. In particular, assigning $\boldsymbol{\gamma}^{\mathrm{LME}}_{\mathrm{IF}} = 0.8$ to repatoms near the interface (within one spacing) and $\boldsymbol{\gamma}^{\mathrm{LME}}_{\mathrm{FF}} = 2.0$ in the far field provides a simple and accurate default for $\boldsymbol{\gamma}^{\mathrm{LME}}$. This choice enables the extended QC method with LME interpolation to capture weak discontinuities in discrete lattices without requiring per-repatom optimization.
\end{itemize}
 
A summation rule consistent with the LME interpolation and Heaviside enrichment remains to be developed, representing a current limitation of the present work. Future studies might extend the approach beyond single-inclusion configurations to more realistic inclusion geometries and heterogeneous discrete lattices with complex microstructures. Moreover, the proposed framework can readily be generalized to three-dimensional settings and to different stiffness contrasts, including cases with soft inclusions embedded in stiff matrices.

%\backmatter
\bmsection*{Author contributions}

{\bf Benjamin Werner}: Methodology, Software, Validation, Formal analysis, Investigation, Writing--Original Draft, Visualization;
{\bf Ond\v{r}ej Roko\v{s}}: Conceptualization, Methodology, Software, Formal analysis, Resources, Writing--Review \& Editing, Supervision, Funding acquisition; {\bf Jan Zeman}: Conceptualization, Methodology, Formal analysis, Writing--Review \& Editing, Project administration, Funding acquisition;

\bmsection*{Acknowledgments}
The work of Benjamin Werner received funding from projects No.~CZ.02.2.69/0.0/0.0/18\_053/\allowbreak 0016980 awarded by the Ministry of Education, Youth and Sports of the Czech Republic (from 02/2021--04/2022), No.~22-35755K awarded the Czech Science Foundation (from 01/2023--08/2024), by the CTU Global Postdoc Fellowship Program (from 05/2022--05/2024), and by the European Union (Marie Sk\l{}odowska-Curie grant agreement No. 101151096, from 09/2024). The work of Jan Zeman and Ond\v{r}ej Roko\v{s} was supported by project No.~19-26143X awarded by the Czech Science Foundation.

\bmsection*{Financial disclosure}

None reported.

\bmsection*{Conflict of interest}

The authors declare no conflicts of interest.

\bibliographystyle{plain}
\bibliography{references_QCMaxEnt}

\bmsection*{Supporting information}

Additional supporting information may be found in the online version of the article at the publisher's website.

\appendix
\small

\bmsection{Lagrange multiplier of LME basis functions}
\label{appendix:app1}
\vspace*{12pt}

For the minimization problem in Eq. (\ref{eq:min_lambda}) the Lagrange multiplier $\boldsymbol{\lambda}$ is determined employing the Newton-Raphson method with the residual $\mathbf{r}^{\mathrm{LME}}$ and the Hessian $\mathbf{J}^{\mathrm{LME}}$ of the objective function $\mathrm{log} \, Z (\mathbf{r}_{0},\boldsymbol{\lambda}, \boldsymbol{\beta}^{\mathrm{LME}})$. For brevity, the arguments are omitted. The new Lagrange multiplier
\begin{equation}
    \boldsymbol{\lambda}_{n+1}
    =
    \boldsymbol{\lambda}_{n} 
    - 
    \left( \mathbf{J}^{\mathrm{LME}} \right) ^{-1}
    \mathbf{r}^{\mathrm{LME}}
\end{equation}
is determined by the value of the previous iteration $\boldsymbol{\lambda}_{n}$, the residual and the Hessian. The residual
\begin{equation}
    \mathbf{r}^{\mathrm{LME}}
    =
    \dfrac{\partial}{\partial \boldsymbol{\lambda}} \mathrm{log} \, Z
    =
    \sum_{\alpha \in N_\mathrm{rep}} 
    \phi_{\alpha} \left(\mathbf{r}_{0} - \mathbf{r}_{\mathrm{0, \, \mathrm{rep}}}^{\alpha}\right)
    \label{eq:LME_residual}
\end{equation}
is the first derivative of the objective function $\mathrm{log} \, Z$ with respect to $\boldsymbol{\lambda}$, while the Hessian 
\begin{equation}
    \mathbf{J}^{\mathrm{LME}}
    =
    \dfrac{\partial^2}{\partial \boldsymbol{\lambda} \partial \boldsymbol{\lambda}} \mathrm{log} \, Z 
    =
    \sum_{\alpha \in N_\mathrm{rep}} 
    \phi_{\alpha} \left(\mathbf{r}_{0} - \mathbf{r}_{\mathrm{0, \, \mathrm{rep}}}^{\alpha}\right)
    \otimes 
    \left(\mathbf{r}_{0} - \mathbf{r}_{\mathrm{0, \, \mathrm{rep}}}^{\alpha}\right)
    -
    \mathbf{r}^{\mathrm{LME}}
    \otimes
    \mathbf{r}^{\mathrm{LME}}
\end{equation}
is the second derivative of the objective function $\mathrm{log} \, Z$. Foca \cite{Foca2015} used the regularized Newton-Raphson method, suggested by Polyak \cite{Polyak2009125}, to solve the minimization problem in Eq. (\ref{eq:min_lambda}) with improved convergence behavior. The objective function is replaced by the strongly convex one
\begin{equation}
    F \left(\mathbf{r}_{0},\boldsymbol{\lambda},\boldsymbol{\zeta}, \boldsymbol{\beta}^{\mathrm{LME}}\right)
    =
    \mathrm{log} \, Z 
    +
    \normtwo{\dfrac{\partial}{\partial \boldsymbol{\lambda}} \mathrm{log} \, Z} \,
    \normtwo{ \boldsymbol{\zeta} - \boldsymbol{\lambda} }^2
\end{equation}
for the regularized Newton-Raphson method, which leads to
\begin{equation}
    \boldsymbol{\lambda}^*(\mathbf{r}_{0})
    =
    \mathrm{arg \, min}_{\zeta \in \mathbb{R}^2} \,F \rvert_{\boldsymbol{\zeta} - \boldsymbol{\lambda}}.
\end{equation}
The first derivative of the new objective function
\begin{equation}
    \dfrac{\partial F}{\partial \boldsymbol{\lambda}} \rvert_{\boldsymbol{\zeta} = \boldsymbol{\lambda}}
    =
    \mathbf{r}^{\mathrm{LME}}
\end{equation}
yields the same residual as in Eq.~(\ref{eq:LME_residual}). The regularized Hessian
\begin{equation}
    \dfrac{\partial^2 F}{\partial \boldsymbol{\lambda} \partial \boldsymbol{\lambda}} \rvert_{\boldsymbol{\zeta} = \boldsymbol{\lambda}}
    =
    \mathbf{J}^{\mathrm{LME}}
    +
    \normtwo{\mathbf{r}^{\mathrm{LME}} } \mathbf{I}
\end{equation} 
has an additional term with the Euclidean norm of the residual multiplied by the identity matrix $\mathbf{I}$. The new Lagrange multiplier using the regularized Newton-Raphson method is therefore calculated by
\begin{equation}
    \boldsymbol{\lambda}_{n+1}
    =
    \boldsymbol{\lambda}_{n} 
    - 
    \left( \mathbf{J}^{\mathrm{LME}} 
    +
    \normtwo{\mathbf{r}^{\mathrm{LME}}} \mathbf{I} \right)^{-1}
    \mathbf{r}^{\mathrm{LME}}.
\end{equation}
This avoids numerical issues due to a singular Hessian $\mathbf{J}^{\mathrm{LME}}$, the solution is independent of the starting values $\boldsymbol{\lambda}_0$, and leads to a faster convergence of the minimization problem.

\bmsection{Derivatives of the total potential energy and the shape functions}
\label{appendix:app2}
\vspace*{12pt}

\bmsubsection{Derivatives of the total potential energy}
\label{app2.1a}

The derivative of the total potential energy with respect to $\beta^{\mathrm{LME}}_{\beta}$ is determined using the chain rule
\begin{equation}
    \dfrac{\partial \Pi}{\partial \beta^{\mathrm{LME}}_{\beta}} 
    = 
    \underbrace{\left. \boldsymbol{\Phi}^T \dfrac{\partial \Pi}{\partial \mathbf{r}} \right|_{{\mathbf{r}} = \boldsymbol{\Phi}\mathbf{r}_{\mathrm{rep}}}}_{\mathbf{G}_{\mathrm{int}}}
    \dfrac{\partial \mathbf{r}}{\partial \beta^{\mathrm{LME}}_{\beta}}
    =
    \left( \dfrac{\partial \boldsymbol{\Phi}}{\partial \beta^{\mathrm{LME}}_{\beta}} \mathbf{r}_{\mathrm{rep}} \right)^T
    \mathbf{G}_{\mathrm{int}}
    \label{eq:dPi_dbeta1}
\end{equation}
where we make use of $\boldsymbol{\Phi}(\boldsymbol{\beta}^{\mathrm{LME}}_{\beta})$ being a function of the locality parameter $\boldsymbol{\beta}^{\mathrm{LME}}_{\beta}$ through individual LME interpolation functions $\phi_{\alpha}$.

\bmsubsection{Derivatives of the interpolation matrix}
\label{app2.1b}

For the derivative of the LME interpolation functions with respect to the locality parameter $\beta^{\mathrm{LME}}_{\alpha}$, the symbol $^*$ is used to indicate a function evaluated at $\boldsymbol{\lambda}^*(\mathbf{r}_{0}, \boldsymbol{\beta}^{\mathrm{LME}}) = \mathrm{arg \, min}_{\zeta \in \mathbb{R}^2} \,F(\mathbf{r}_{0},\boldsymbol{\lambda}, \boldsymbol{\zeta}, \boldsymbol{\beta}) \rvert_{\boldsymbol{\zeta} - \boldsymbol{\lambda}}$, i.e., at the optimum. The derivative of the LME interpolation function was presented by Rosolen et al. \cite{Rosolen2010868} and is stated here for completeness. The basis functions can be written as
\begin{align}
    \phi_{\alpha}(\mathbf{r}_{0},\lambda,\boldsymbol{\beta}^{\mathrm{LME}}) 
    =
    \dfrac{\mathrm{exp}[f_{\alpha}(\mathbf{r}_{0},\lambda,\beta^{\mathrm{LME}}_{\alpha}) ]}{\sum_{\beta}\mathrm{exp}[f_{\beta}(\mathbf{r}_{0},\lambda,\beta^{\mathrm{LME}}_{\beta})]}
    =
    \dfrac{\mathrm{exp}[f_{\alpha}(\mathbf{r}_{0},\lambda,\beta^{\mathrm{LME}}_{\alpha}) ]}{Z(\mathbf{r}_{0},\lambda,\boldsymbol{\beta}^{\mathrm{LME}})}
\end{align}
with
\begin{align}
    f_{\alpha}(\mathbf{r}_{0},\lambda,\beta^{\mathrm{LME}}_{\alpha}) 
    =
    - \beta^{\mathrm{LME}}_{\alpha} \normtwo{\mathbf{r}_{0} - \mathbf{r}_{\mathrm{0, \, rep}}^{\alpha}}^2
    + \lambda \left(\mathbf{r}_{0} - \mathbf{r}_{\mathrm{0, \, rep}}^{\alpha}\right)
\end{align}
The derivative
\begin{align}
% Eq 1
    \dfrac{\partial \phi_{\alpha}^*}{\partial \beta^{\mathrm{LME}}_{\beta}} 
    = 
    \phi_{\alpha}^*\left(\dfrac{\partial f_{\alpha}^*}{\partial \beta^{\mathrm{LME}}_{\beta}} - \sum_{\gamma} \phi_{\gamma}^*\dfrac{\partial f_{\gamma}^*}{\partial \beta^{\mathrm{LME}}_{\beta}}\right)
    \label{eq:dpa*_dbetab}
\end{align}
with
\begin{align}
% Eq 2
    \dfrac{\partial f_{\alpha}^*}{\partial \beta_{\beta}^{\mathrm{LME}}} 
    = 
    \left(\dfrac{\partial f_{\alpha}}{\partial \beta_{\beta}^{\mathrm{LME}}} \right)^* + \left(\dfrac{\partial f_{\alpha}}{\partial \lambda} \right)^*\left(\dfrac{\partial \lambda}{\partial \beta_{\beta}^{\mathrm{LME}}} \right)^*
    \label{eq:dfa*_dbetab}
\end{align}
where the first term in Eq. (\ref{eq:dfa*_dbetab}) is explicitly expressed as
\begin{align}
% Eq 3
     \left(\dfrac{\partial f_{\alpha}}{\partial \beta_{\beta}^{\mathrm{LME}}} \right)^* 
     = 
     -\delta_{\alpha \beta} \normtwo{\mathbf{r}_{0} - \mathbf{r}_{\mathrm{0, \, rep}}^{\alpha}}^2,
     \label{eq:dfa_dbeta_star}
\end{align}
and the second term is comprised by
\begin{align}
% Eq 4
     \left(\dfrac{\partial f_{\alpha}}{\partial \lambda} \right)^*
     = 
     \left(\mathbf{r}_{0} - \mathbf{r}_{\mathrm{0, \, rep}}^{\alpha}\right)
     \label{eq:dfa_dlambda_star}
\end{align}
and 
\begin{align}
% Eq 5
    \left(\dfrac{\partial \lambda}{\partial \beta_{\beta}^{\mathrm{LME}}} \right)^* 
    = 
    \phi_{\beta}^* \normtwo{\mathbf{r}_{0} - \mathbf{r}_{\mathrm{0, \, rep}}^{\beta}}^2
    \left( \mathbf{J}^{\mathrm{LME}\,*}\right)^{-1}
    \left(\mathbf{r}_{0} - \mathbf{r}_{\mathrm{0, \, rep}}^{\beta}\right).
    \label{eq:dlambda_dbeta_star}
\end{align}
Here, $\mathbf{J}^{\mathrm{LME}\,*}$ denotes the Hessian $\mathbf{J}^{\mathrm{LME}}$ evaluated at $\boldsymbol{\lambda}^*(\mathbf{r}_0,\boldsymbol{\beta}^{\mathrm{LME}})$. Substituting these results back into Eq. (\ref{eq:dfa*_dbetab}) leads to
\begin{align}
% Eq 6
    \dfrac{\partial f_{\alpha}^*}{\partial \beta_{\beta}^{\mathrm{LME}}} 
    = 
    -\delta_{\alpha \beta}\normtwo{\mathbf{r}_{0} - \mathbf{r}_{\mathrm{0, \, rep}}^{\alpha}}^2
    +
    \phi_{\beta}^*\normtwo{\mathbf{r}_{0} - \mathbf{r}_{\mathrm{0 \, rep}}^{\beta}}^2 \left(\mathbf{r}_{0} - \mathbf{r}_{\mathrm{0, \, rep}}^{\alpha}\right) \cdot 
    \left( \mathbf{J}^{\mathrm{LME}\,*}\right)^{-1}
    \left(\mathbf{r}_{0} - \mathbf{r}_{\mathrm{0, \, rep}}^{\beta}\right).
    \label{eq:fa_beta}
\end{align}
The same equation can be used for the second term in Eq. (\ref{eq:dpa*_dbetab}) simply replacing index $\alpha$ by $\gamma$, leading to
\begin{align}
% Eq 7
    \dfrac{\partial f_{\gamma}^*}{\partial \beta_{\beta}^{\mathrm{LME}}} 
    = 
    -\delta_{\gamma \beta}\normtwo{\mathbf{r}_{0} - \mathbf{r}_{\mathrm{0, \, rep}}^{\gamma}}^2
    +
    \phi_{\beta}^*\normtwo{\mathbf{r}_{0} - \mathbf{r}_{\mathrm{0, \, rep}}^{\beta}}^2 \left(\mathbf{r}_{0} - \mathbf{r}_{\mathrm{0, \, rep}}^{\gamma}\right) \cdot 
    \left( \mathbf{J}^{\mathrm{LME}\,*}\right)^{-1}
    \left(\mathbf{r}_{0} - \mathbf{r}_{\mathrm{0, \, rep}}^{\beta}\right).
\end{align}
Summing over $\gamma$ and multiplying with $\phi_{\gamma}^*$ gives
\begin{align}
% Eq 8
    \sum_{\gamma} \phi_{\gamma}^*\dfrac{\partial f_{\gamma}^*}{\partial \beta_{\beta}^{\mathrm{LME}}}
    = 
    \sum_{\gamma} \phi_{\gamma}^* \left(
    -\delta_{\gamma \beta} \normtwo{\mathbf{r}_{0} - \mathbf{r}_{\mathrm{0, \, rep}}^{\gamma}}^2
    +
    \phi_{\beta}^* \normtwo{\mathbf{r}_{0} - \mathbf{r}_{\mathrm{0, \, rep}}^{\beta}}^2 \left(\mathbf{r}_{0} - \mathbf{r}_{\mathrm{0 \, rep}}^{\gamma}\right) \cdot 
    \left( \mathbf{J}^{\mathrm{LME}\,*}\right)^{-1}
    \left(\mathbf{r}_{0} - \mathbf{r}_{\mathrm{0, \, rep}}^{\beta}\right)\right),
\end{align}
which can be simplified to
\begin{align}
% Eq 9
    \sum_{\gamma} \phi_{\gamma}^*\dfrac{\partial f_{\gamma}^*}{\partial \beta_{\beta}^{\mathrm{LME}}}
    = 
    - \phi_{\beta}^* \normtwo{\mathbf{r}_{0} - \mathbf{r}_{\mathrm{0, \, rep}}^{\beta}}^2
    +
    \underbrace{
    \sum_{\gamma} \phi_{\gamma}^* \left(
    \phi_{\beta}^* \normtwo{\mathbf{r}_{0} - \mathbf{r}_{\mathrm{0, \, rep}}^{\beta}}^2 
    \left(\mathbf{r}_{0} - \mathbf{r}_{\mathrm{0, \, rep}}^{\gamma}\right) 
    \cdot 
    \left( \mathbf{J}^{\mathrm{LME}\,*}\right)^{-1}
    \left(\mathbf{r}_{0} - \mathbf{r}_{\mathrm{0, \, rep}}^{\beta}\right) 
    \right)
    }_\text{$\approx 0$ and neglected}.
    \label{eq:pc_fc_beta}
\end{align}
Substituting Eqs.~(\ref{eq:pc_fc_beta}) and (\ref{eq:fa_beta}) into Eq.~(\ref{eq:dpa*_dbetab}) yields
\begin{align}
% Eq 10
    \dfrac{\partial \phi_{\alpha}^*}{\partial \beta_{\beta}^{\mathrm{LME}}} 
    = 
    \phi_{\alpha}^*\left(
    -\delta_{\alpha \beta} \normtwo{\mathbf{r}_{0} - \mathbf{r}_{\mathrm{0, \, rep}}^{\alpha}}^2
    +
    \phi_{\beta}^* \normtwo{\mathbf{r}_{0} - \mathbf{r}_{\mathrm{0, \, rep}}^{\beta}}^2 
    \left(\mathbf{r}_{0} - \mathbf{r}_{\mathrm{0, \, rep}}^{\alpha}\right) 
    \cdot 
    \left( \mathbf{J}^{\mathrm{LME}\,*}\right)^{-1}
    \left(\mathbf{r}_{0} - \mathbf{r}_{\mathrm{0, \, rep}}^{\beta}\right)
    + 
    \phi_{\beta}^* \normtwo{\mathbf{r}_{0} - \mathbf{r}_{\mathrm{0, \, rep}}^{\beta}}^2
    \right).
    \label{eq:dpa*_dbetab_2}
\end{align}
which simplifies to
\begin{align}
% Eq 11
    \dfrac{\partial \phi_{\alpha}^*}{\partial \beta^{\mathrm{LME}}_{\beta}} 
    = 
    \phi_{\alpha}^* \normtwo{\mathbf{r}_{0} - \mathbf{r}_{\mathrm{0, \, rep}}^{\beta}}^2 
    \left(
    \phi_{\beta}^* \left[ \left(\mathbf{r}_{0} - \mathbf{r}_{\mathrm{0, \, rep}}^{\alpha}\right) 
    \cdot 
    \left( \mathbf{J}^{\mathrm{LME}\,*}\right)^{-1}
    \left(\mathbf{r}_{0} - \mathbf{r}_{\mathrm{0, \, rep}}^{\beta}\right)
    + 
    1\right]
    -
    \delta_{\alpha \beta}
    \right).
    \label{eq:dpa_star_dbeta}
\end{align}

\bmsubsection{Derivatives of the enriched interpolation matrix}
\label{app2.1c}

The derivative of the enriched interpolation functions $\phi_{ \mathrm{GSON}, \, j}^{\star}$ with respect to $\beta^{\mathrm{LME}}_{\beta}$ is calculated by the derivative of Eq. (\ref{eq:GS_basis_func})
\begin{equation}
    \dfrac{\partial \phi_{ \mathrm{GS}, \, j}^{\star}}{\partial \beta^{\mathrm{LME}}_{\beta_j}} 
    = 
    \dfrac{\partial \phi_j^{\star}}{\partial \beta^{\mathrm{LME}}_{\beta_j}} 
    - 
    \sum_{i=1}^{j-1} 
    \left[
    \left(
    \dfrac{\partial \phi_j^{\star}}{\partial \beta^{\mathrm{LME}}_{\beta_j}}^T \cdot \phi_{\mathrm{GSON}, \, i}^{\star}
    +
    \phi_j^{\star \, T} \cdot \dfrac{\partial \phi_{\mathrm{GSON}, \, i}^{\star}}{\partial \beta^{\mathrm{LME}}_{\beta_i}}
    \right) \phi_{\mathrm{GSON}, \, i}^{\star}
    +
    \left( \phi_j^{\star \, T} \cdot \phi_{\mathrm{GSON}, \, i}^{\star} \right) \dfrac{\partial \phi_{\mathrm{GSON}, \, i}^{\star}}{\partial \beta^{\mathrm{LME}}_{\beta_i}}
    \right]
\end{equation}
where the first term is
\begin{equation}
    \dfrac{\partial \phi_{j}^{\star}}{\partial \beta^{\mathrm{LME}}_{\beta_j}}
    = 
    \dfrac{\partial \phi_{j}}{\partial \beta^{\mathrm{LME}}_{\beta_j}} 
    \left(\chi(\mathbf{r}_0^{\alpha}) - \chi(\mathbf{r}_0^{\beta_j}) \right)
\end{equation}
with the shifted Heaviside enrichment $\chi(\mathbf{r}_0^{\alpha}) - \chi(\mathbf{r}_0^{\beta_j})$. The derivative of the orthonormalized Gram--Schmidt interpolation functions
\begin{equation}
    \dfrac{\partial \phi_{\mathrm{GSON}, \, j}^{\star} }{\partial \beta^{\mathrm{LME}}_{\beta_j}} 
    =
    \dfrac{\dfrac{\partial \phi_{\mathrm{GS}, \, j}^{\star}}{\partial \beta^{\mathrm{LME}}_{\beta_j}} ||\phi_{\mathrm{GS}, \, j}^{\star}|| 
    - 
    \dfrac{\partial ||\phi_{\mathrm{GS}, \, j}^{\star}||}{\partial \beta^{\mathrm{LME}}_{\beta_j}}   \phi_{\mathrm{GS}, \, j}^{\star}}{||\phi_{\mathrm{GS}, \, j}^{\star}||^2}
\end{equation}
is determined from Eq. (\ref{eq:GSON_basis_func}) using the expression
\begin{equation}
    \dfrac{\partial}{\partial \beta^{\mathrm{LME}}_{\beta_j}} ||\phi_{\mathrm{GS}, \, j}^{\star}|| 
    = 
    \dfrac{\phi_{\mathrm{GS}, \, j}^{\star \, T} \cdot \dfrac{\partial \phi_{\mathrm{GS}, \, j}^{\star}}{\partial \beta^{\mathrm{LME}}_{\beta_j}} }{||\phi_{\mathrm{GS}, \, j}^{\star}||}.
\end{equation}

\end{document}